\documentclass[11pt]{amsart}
\usepackage{amsmath,amsfonts,amsthm,amssymb,url,graphicx,algorithm, algorithmicx,tikz,sagetex,float,hyperref}
\usepackage[noend]{algpseudocode}
\usetikzlibrary{arrows,automata,positioning}

\DeclareFontFamily{OT1}{rsfs}{}
\DeclareFontShape{OT1}{rsfs}{n}{it}{<-> rsfs10}{}
\DeclareMathAlphabet{\mathscr}{OT1}{rsfs}{n}{it}
\DeclareMathOperator{\sgn}{sgn}

\DeclareMathOperator{\len}{len}

\DeclareMathOperator{\den}{den}

\DeclareMathOperator{\mo}{\,mod}

\newtheorem{prop}{Proposition}[section]

\newtheorem*{main}{Main Theorem}
 
\newtheorem{lem}[prop]{Lemma}

\newtheorem*{defn*}{Definition}

\numberwithin{equation}{section}
\title{Summing $\mu(n)$: a faster elementary algorithm}
\author{Harald Andr\'es Helfgott}
\address{Harald A. Helfgott, 
  Mathematisches Institut,
Georg-August Universit\"{a}t G\"{o}ttingen, Bunsenstra{\ss}e 3-5, D-37073 G\"{o}ttingen,
Germany; IMJ-PRG, UMR 7586,
  58 avenue de France, B\^{a}timent S. Germain, case 7012,
  75013 Paris CEDEX 13, France}
\email{harald.helfgott@gmail.com}
\author{Lola Thompson}
\address{Lola Thompson, 
Mathematics Institute, Utrecht University,  Hans Freudenthalgebouw, Budapestlaan 6, 3584 CD Utrecht, Netherlands}
\email{l.thompson@uu.nl}
\begin{document}
\begin{abstract}


We present a new elementary algorithm that takes 
\[\begin{aligned}
\mathrm{time} \ \ O_\epsilon\left(x^{\frac{3}{5}} (\log x)^{\frac{3}{5}+\epsilon} \right)
\ \ \mathrm{and}\ \ \mathrm{space} \ \ O\left(x^{\frac{3}{10}} (\log x)^{\frac{13}{10}}
\right)\end{aligned}\]
for computing $M(x) = \sum_{n \leq x} \mu(n),$ where $\mu(n)$ is the M\"{o}bius function. This is the first improvement in the exponent of $x$ for an elementary algorithm since 1985.

We also show that it is possible to reduce
space consumption to $O(x^{1/5} (\log x)^{5/3})$ by the use of (Helfgott, 2020), at the cost of letting time
rise to the order of $x^{3/5} (\log x)$.
\end{abstract}

\maketitle

\section{Introduction}

There are several well-studied sums in analytic number theory that involve the M\"{o}bius function. For example, Mertens \cite{Mertens}
considered $$M(x) = \sum_{n \leq x} \mu(n),$$
now called the {\em Mertens function}. Based on numerical evidence, he conjectured that $|M(x)| \leq \sqrt{x}$ for all $x$. His conjecture was
disproved by Odlyzko and te Riele \cite{OtR1985}. Pintz \cite{Pintz87} made their result effective, showing that there exists a value of $x < \exp(3.21 \times 10^{64})$ for which $|M(x)| > \sqrt{x}$. It is still not known when $|M(x)| > \sqrt{x}$ holds for the first time; Dress \cite{Dress1993} has shown that it cannot hold for $x\leq 10^{12}$, and Hurst has carried out a verification up to $10^{16}$  \cite{Hurst18}. Isolated values of $M(x)$ have been computed in \cite{Dress1993} and in subsequent papers.  

The two most time-efficient algorithms known for computing $M(x)$ are the following:
\begin{enumerate}
\item An analytic algorithm (Lagarias-Odlyzko \cite{LO87}), with computations based
  on integrals of $\zeta(s)$; its running time is $O(x^{1/2+\epsilon})$.
\item A more elementary algorithm (Meissel-Lehmer \cite{Lehmer59} and
  Lagarias-Miller-Odlyzko \cite{LMO}; refined by Del\'eglise-Rivat \cite{DelegliseRivat96}),
  with running time about $O(x^{2/3})$. 
\end{enumerate}
These algorithms are variants of similar algorithms for computing $\pi(x)$, the number of primes up to $x$. The analytic algorithm had to wait for almost 30 years to receive its first
rigorous, unconditional implementation due to Platt \cite{Plattpi}, which concerns only the
computation of $\pi(x)$. The computation of $M(x)$ using the analytic algorithm presents additional complications and has not been implemented. Moreover,
in the range explored to date ($x\leq 10^{22}$), elementary algorithms are faster in practice, at least for computing $\pi(x)$.

Del\'{e}glise and Rivat's paper \cite{DelegliseRivat96} gives the
values of $M(x)$ for $x=10^6,10^7,\dotsc,10^{16}$.
An unpublished 2011 preprint of Kuznetsov \cite{Kuznetsov11} gives the values of $M(x)$ for
$x=10^{16},10^{17},\dotsc,10^{22}$ using parallel computing.
More recently, Hurst \cite{Hurst18} computed $M(x)$  for
 $x = 2^n$, $n\leq 73$. (Note that $2^{73} = 9.444\dotsc \cdot 10^{21}$.)
The computations in \cite{Kuznetsov11} and \cite{Hurst18} are both based on the algorithm in \cite{DelegliseRivat96}. 

Since 1996, all work on these problems has centered on improving the implementation, with no essential improvements to the algorithm or to its computational complexity. The goal of the present paper is to develop a new elementary algorithm that is more time-efficient and space-efficient than the algorithm in \cite{DelegliseRivat96}. We show:

\begin{main} We can compute $M(x)$ in
  \[\mathrm{time} \ \ O\left(x^{\frac{3}{5}} (\log x)^{\frac{3}{5}}
  (\log \log x)^{\frac{2}{5}} \right) \ \
  \mathrm{and \ space} \ \ O\left(x^{\frac{3}{10}} (\log x)^{\frac{13}{10}}
  (\log \log x)^{-\frac{3}{10}} \right).\]
\end{main}

This is the first improvement in the exponent of $x$ since 1985. Using our algorithm, we have been able to extend the work of Hurst and Kuznetsov, computing $M(x)$ for $x = 2^n$, $n\leq 75$, and for $x = 10^{n}$, $n\leq 23$. We expect that professional programmers who have access to significant computer resources will be able to extend this range further. 

\subsection{Our approach}
The general idea used in all of the elementary algorithms
(\cite{LMO}, \cite{DelegliseRivat96}, etc.) is as follows. One always starts with a combinatorial identity to break $M(x)$ into smaller sums.
For example, a variant of Vaughan's identity allows one to rewrite $M(x)$
as follows:
$$M(x) = 2M(\sqrt{x}) - \sum_{n \leq x} \sum_{\substack{m_1 m_2 n_1 = n \\ m_1, m_2 \leq \sqrt{x}}} \mu(m_1) \mu(m_2).$$ Swapping the order of summation, one can write $$M(x) = 2M(\sqrt{x}) - \sum_{m_1, m_2 \leq \sqrt{x}} \mu(m_1) \mu(m_2) \left\lfloor \frac{x}{m_1 m_2} \right\rfloor.$$ The first term can be easily computed in time $O(\sqrt{x} \log \log x)$ and space $O(x^{1/4})$, or else, proceeding as in \cite{Helfgott2020}, in time $O(\sqrt{x} \log x)$ and space $O(x^{1/6} (\log x)^{2/3})$. To handle the subtracted term, the idea is to fix a parameter $v \leq \sqrt{x}$, and then split the sum into two sums: one over $m_1, m_2 \leq v$ and the other with $\max(m_1, m_2) > v$. The difference between the approach taken in the present paper and those that came before it is that our predecessors take $v = x^{1/3}$ and then compute the sum for $m_1, m_2 \leq v$ in time $O(v^2)$. We will take our $v$ to be a
little larger, namely, about $x^{2/5}$. Because we take a larger value of $v$, we have
to treat the case with $m_1, m_2 \leq v$ with greater care than
\cite{DelegliseRivat96} et al. Indeed, the bulk of our work will be in Section
\ref{sec:largefree}, where we show how to handle this case. 

Our approach in Section \ref{sec:largefree} roughly amounts to analyzing the difference between reality and
a model that we obtain via Diophantine approximation, in that we show
that this difference
has a simple description in terms of congruence classes and segments.
This description allows us to compute the difference quickly, in part
by means of table lookups. 


\subsection{Alternatives} In a previous draft of our paper, we followed a route
more closely related to the main ideas in papers by Galway \cite{Galway} and by the first author \cite{Helfgott2020}. Those papers
succeeded in reducing the space needed for implementing the sieve of Eratosthenes (or the Atkin-Bernstein sieve, in Galway's case) down to about $O(x^{1/3})$.
In particular, \cite{Helfgott2020} provides an algorithm for computing
$\mu(n)$ for all successive $n\leq x$ in time $O(x \log x)$ and space $O(x^{1/3} (\log x)^{2/3})$, building on an approach from a paper of Croot, Helfgott, and Tao \cite{TCH12} that computes $\sum_{n\leq x} \tau(n)$ in time about $O(x^{1/3})$. That approach is in turn related to Vinogradov's take on the divisor problem \cite[Ch.~III, exer.~3-6]{Vinogradov}
(based on Vorono\"{i}).

The total time taken by the algorithm in the previous version of
our paper was on the order of $x^{3/5} (\log x)^{8/5}$.
Thus, the current version is asymptotically faster.
If an unrelated improvement present
in the current version (Algorithm \ref{alg:factofun};
see \S \ref{sec:largenonfree}) were introduced in the older version,
time usage would be on the order of $x^{3/5} (\log x)^{6/5}
(\log \log x)^{2/5}$. 
We sketch the older version of the algorithm in Appendix 
\ref{subs:appalt}.

Of course, we could use \cite{Helfgott2020} as a black box to reduce space
consumption in some of our routines, while leaving everything else as it
is in the current version. Time complexity would increase slightly,
while space complexity would
be much reduced.
 More precisely: using \cite{Helfgott2020} as a black box, and keeping
everything else the same,
we could compute $M(x)$  in time $O(x^{3/5} (\log x))$
and space $O(x^{1/5} (\log x)^{5/3})$. We choose to focus instead
on the version of the algorithm reflected in the main theorem;
it is faster but less space-efficient.




\subsection{Notation and algorithmic conventions} As usual, we write $f(x) = O(g(x))$ to denote that there is a positive constant $C$ such that $|f(x)| \leq C g(x)$ for all sufficiently large $x$. The notation $f(x)\ll g(x)$ is synonymous
to $f(x)=O(g(x))$.
We use $f(x) = O^*(g(x))$ to indicate  something stronger, namely, $|f(x)| \leq g(x)$ for all $x$.

For $x \in \mathbb{R}$, we write $\lfloor x \rfloor$ for the largest integer $\leq x$, and $\{x\}$ for $x - \lfloor x \rfloor$. Thus, $\{x\} \in [0, 1)$ no matter whether $x < 0$, $x = 0$, or $x > 0$.

  We write $\log_b x$ to mean the logarithm base $b$ of $x$, {\em not}
  $\log \log \cdots \log x$ ($\log$ iterated $b$ times).
  
  Throughout this paper, we assume that arithmetic operations take time $O(1)$, and we count space in bits. The combination of these two assumptions may seem
  counterintuitive, but it is actually a good reflection of practice,
  particularly given that any $x$ for which we can compute $M(x)$ in
  reasonable time can be stored in a fixed-sized integer (64 or 128 bits). All of the pseudocode for our algorithms appears at the end of this paper.


  \subsection{Acknowledgements} The authors would like to thank the Max Planck Institute for Mathematics, which hosted the two of them for a joint visit from February 1 - April 15, 2020. They are especially grateful to have had access to the parallel computers at the MPIM. While completing this research,
H.\nobreak\hspace{0.1em}H.\
was partially supported by the European Research Council under Programme H2020-EU.1.1., ERC Grant ID: 648329 (codename GRANT), and by his Humboldt professorship. L.\nobreak\hspace{0.1em}T.\
was partially supported by the Max Planck Institute for Mathematics for her sabbatical during the 2019 - 2020 academic year. This work began while she was employed by Oberlin College. She is grateful to Oberlin for supporting her during the early stages of this project.

\section{Preparatory work: identities}

We will start from the identity
\begin{equation}\label{eq:firstid}
  \mu(n) = - \mathop{\sum_{m_1 m_2 n_1 = n}}_{m_1,m_2\leq u} \mu(m_1) \mu(m_2)
+ \begin{cases} 2 \mu(n) &\text{if $n\leq u$}\\ 0 &\text{otherwise,}
  \end{cases}\end{equation}
valid for $n\leq x$ and $u\geq \sqrt{x}$. (We will set $u=\sqrt{x}$.)
This identity is simply the case $K=2$ of
Heath-Brown's identity 
for the M\"obius function: for all $K \geq 1, n \geq 1$, and $u \geq n^{1/K}$,
$$\mu(n) = - \sum_{1 \leq k \leq K} (-1)^k \binom{K}{k}
\mathop{\sum_{m_1 ... m_k n_1 ... n_{k-1}=n}}_{m_1,...,m_k\leq u} \mu(m_1) ... \mu(m_k).$$
(See \cite[(13.38)]{MR2061214}; note, however, that there is a typographical error
under the sum there: $m_1\dotsc m_k n_1\dotsc n_{k}=n$ should be
$m_1\dotsc m_k n_1\dotsc n_{k-1}=n$.)
Alternatively, we can derive (\ref{eq:firstid}) immediately from
Vaughan's identity for $\mu$: that identity would, in general, have
a term consisting of a sum over all decompositions $m_1 m_2 n_1 = n$ with
$m_1,m_2>u$, but
that term is empty because $u^2\geq x$.

We sum over all $n \leq x$, and obtain
\begin{equation}\label{eq:baseid}
  \begin{aligned} M(x) = 2 M(u) - \sum_{n\leq x} \mathop{\sum_{m_1 m_2 n_1 = n}}_{m_1, m_2\leq u}
  \mu(m_1) \mu(m_2) .
\end{aligned}\end{equation} for $u\geq \sqrt{x}$.

Before we proceed, let us compare matters to the initial approach
 in \cite{DelegliseRivat96}. Lemma 2.1 in 
\cite{DelegliseRivat96} states that
\begin{align}\label{eq:lehman} M(x) = M(u) - \sum_{m \leq u} \mu(m) \sum_{\frac{u}{m} < n \leq \frac{x}{m}} M\left(\frac{x}{mn}\right)\end{align} for
$1\leq u\leq x$. This identity is due to
Lehman \cite[p. 314]{lehman1960liouville}; like Vaughan's identity, it can be proved
essentially by M\"obius inversion.
For $u=\sqrt{x}$, this identity is equivalent to (\ref{eq:firstid}),
as we can see by a change of variables and, again, M\"obius inversion.

We will set $u=\sqrt{x}$ once and for all. We can compute $M(u)$ in
(\ref{eq:baseid}) in time $O(u \log \log u)$ and space $O(\sqrt{u})$,
by a segmented sieve of Eratosthenes. (Alternatively,
we can compute $M(u)$
in time $O(u \log u)$ and space $O(u^{1/3} (\log u)^{2/3})$, using
the space-optimized version of the segmented sieve of
Eratosthenes in \cite{Helfgott2020}.) Thus, we will be able to focus
on the other term on the right side of (\ref{eq:baseid}).
We can write, for any $v\leq u$,
\begin{equation}\label{eq:nemero}\begin{aligned}
 \sum_{n\leq x} \mathop{\sum_{m_1 m_2 n_1 = n}}_{m_1, m_2\leq u} \mu(m_1) \mu(m_2) =
 &\sum_{n\leq x}\; \mathop{\sum_{m_1 m_2 n_1 = n}}_{m_1, m_2\leq v}  \mu(m_1) \mu(m_2)\\
 + &\sum_{n\leq x} \mathop{\mathop{\sum_{m_1 m_2 n_1 = n}}_{m_1, m_2\leq u}}_{\max(
   m_1,m_2)>v} \mu(m_1) \mu(m_2).
\end{aligned}\end{equation}

In this way, computing $M(x)$ reduces to computing the two double sums
on the right side of (\ref{eq:nemero}).

\section{The case of a large non-free variable}\label{sec:largenonfree}

Let us work on the second sum in (\ref{eq:nemero}) first.
It is not particularly difficult to deal with;
there are a few alternative procedures
that would lead to the same time complexity,
and several that would lead to a treatment whose time complexity is
worse by only a factor of $\log x$.

Clearly,
\begin{equation}\label{eq:dullea}\begin{aligned}
\sum_{n\leq x} \mathop{\mathop{\sum_{m_1 m_2 n_1 = n}}_{m_1, m_2\leq u}}_{\max(
  m_1,m_2)>v} \mu(m_1) \mu(m_2) &=
\sum_{v<m\leq u} \mu(m)^2 \left\lfloor \frac{x}{m^2}\right\rfloor \\ &+
2 \sum_{n\leq x} \mathop{\mathop{\sum_{m_1 m_2 n_1 = n}}_{v < m_1\leq u}}_{m_2<m_1}
\mu(m_1) \mu(m_2)
\end{aligned}\end{equation}
and
\begin{equation}\label{eq:pina}
  \sum_{n\leq x} \mathop{\mathop{\sum_{m_1 m_2 n_1 = n}}_{v < m_1\leq u}}_{m_2<m_1}
\mu(m_1) \mu(m_2) =
\sum_{v<a\leq u} \mu(a) \sum_{r\leq \frac{x}{a}}
\mathop{\sum_{b|r}}_{b<a} \mu(b).\end{equation}
It is evident that the first sum on the right in (\ref{eq:dullea})
can be computed in time $O(u \log \log u)$ and space $O(\sqrt{u})$, again by a segmented
sieve. (Alternatively, we can compute it in time $O(u \log u)$ and space $O(u^{1/3} (\log u)^{2/3})$, using the segmented sieve in \cite{Helfgott2020}.)

Write $D(r,y) =
\sum_{b|r: b\leq y} \mu(b)$. Then
\[\begin{aligned}
\sum_{r\leq \frac{x}{a}} \mathop{\sum_{b|r}}_{b<a} \mu(b)
&= \sum_{r\leq \frac{x}{a}} \mathop{\sum_{b|r}}_{b\leq \frac{x}{r}} \mu(b)
- \sum_{r\leq \frac{x}{a}} \mathop{\sum_{b|r}}_{a\leq b\leq \frac{x}{r}} \mu(b)\\
&= \sum_{r\leq \frac{x}{a}} D\left(r,\frac{x}{r}\right)
- \sum_{b\geq a} \mu(b) \sum_{r\leq \frac{x}{b}} 1 =
S\left(\frac{x}{a}\right)
 - \sum_{b\geq a} \mu(b) \left\lfloor \frac{x}{b^2}\right\rfloor.\end{aligned}\]
where 
$S(m) = \sum_{r\leq m} D(r;x/r) = 1 + \sum_{x/u<r\leq m} D(r;x/r)$, since
$D(r;x/r) = \sum_{b|r: b\leq x/r} \mu(b) = \sum_{b|r} \mu(r)$ for
$r\leq \sqrt{x} = u$.

Thus, to compute the right side of (\ref{eq:pina}), it makes sense to let
$n$ take the values
$\lfloor u\rfloor, \lfloor u\rfloor-1,\dotsc, \lfloor v\rfloor+1$ in
descending order; as $n$ decreases, $x/n$ increases, and we compute
$D(r;x/r)$, and thus $S(x/n)$, for increasing values of $r$.
Computing all values of $\mu(a)$ for $v<a\leq u$
using a segmented sieve of Eratosthenes takes time $O(u\log \log u)$
and space $O(\sqrt{u})$.


The main question is how to compute $D(r;x/r)$ efficiently for all $r$
in a given segment.
Using a segmented sieve of Eratosthenes, we can determine the set of
prime
divisors of all $r$ in an interval of the form $\lbrack y, y+\Delta\rbrack$,
$|\Delta|\geq \sqrt{y}$, in time $O(\Delta \log \log y)$ and space $O(\Delta \log y)$.
We want to compute the sum $D(r;x/r) = \sum_{b|r: b<x/r} \mu(b)$ for all
$r$ in that interval. The naive approach would be to
go over all divisors $b$ of all integers $r$ in $\lbrack y, y+\Delta\rbrack$;
since those integers have $\log y$ divisors on average, doing so would take
time $O(\Delta \log y)$.
Fortunately, there is a less obvious way to compute $D(r;x/r)$
in average time $O(\log \log y)$.
We will need a simple lemma on the anatomy of integers.

\begin{lem}\label{lem:isis} Let $P_z(n) = \prod_{p \leq z: p \mid n} p.$ For $z, N, a$ arbitrary and $N < n \leq 2N$ random, the expected value of
  \begin{equation}\label{eq:preteri}
    \mathop{\sum_{\frac{a}{P_z(n)} < d \leq 2a}}_{p \mid d \Rightarrow p > z} \hspace{0.2 in}
  \sum_{\substack{d' \mid n:\; \text{$d'$ squarefree} \\ p \mid d' \Rightarrow z^{1/2} < p \leq z}} 1 
  \end{equation}
is $O(1)$. \end{lem}

\begin{proof} 

  For any fixed positive integer $K$, the numbers $N<n \leq 2 N$ with $P_z(n) = K$ are of the form $m \cdot \prod_{p \leq z: p \mid n}  = m \cdot K,$ where $m$ can be any of the $z$-rough integers $N/K <m \leq  2 N/K$.
  Let us consider how many divisors $d|m$ with properties
  with $p \mid d \Rightarrow p > z$ and $\frac{a}{P_z(n)} < d \leq 2a$
there are on average as $m$
  varies on $(N/K,2 N/K]$.
  
    We can assume that $z\leq N/K$, as otherwise $m$ has at most
    $2$ divisors $d$ free of prime factors $\leq z$ (namely, $d=1$ and $d=m$).
  Then a random integer $m \in (N/K, 2N/K]$ with no prime factors $\leq z$ has the following expected number of divisors in $(\frac{a}{K}, 2a]$:
      $$\frac{1}{(N/K)/\log z} O\left(\sum_{\substack{\frac{a}{K} < d \leq 2a \\ p \mid d \Rightarrow p > z}} \frac{(N/K)/d}{\log z}\right) + O(1)
      = O\Big(1 + 
 \sum_{\substack{\frac{a}{K} < d \leq 2a \\ p \mid d \Rightarrow p > z}} \frac{1}{d}\Big)
,$$
      since the number of integers in $(M, 2M]$ with no prime factors up to $z$ is $\gg M/\log z$ for $z\leq M$ and $\ll M/\log z$ for $z>1$ and $M\geq 1$.
        (The term $O(1)$ is there to account for
        $d=m$; in that case and only then, $(N/K)/d < 1$.)
        

        Applying an upper bound sieve followed by partial summation, we see that
        $$\sum_{\substack{\frac{a}{K} < d \leq 2a \\ p \mid d \Rightarrow p> z}} \frac{1}{d} \ll (\log 2a - \log a/K) \prod_{p \leq z} \left(1 - \frac{1}{p}\right) + 1.$$
        (The term $O(1)$ comes from $\sum_{a/K<d\leq z a/K} 1/d$.)
        By Mertens' Theorem, the product is $\ll 1/\log z$. Hence, $$\sum_{\substack{\frac{a}{K} < d \leq 2a \\ e \mid d \Rightarrow e > z}} \frac{1}{d} = O\left(
        \frac{\log 2 a - \log a/K}{\log z} + 1\right) = O\left(\frac{\log 2 K}{\log z}
         + 1\right).$$

         The number of divisors $d'|n$ with
         $p|d'\Rightarrow z^{1/2}<p\leq z$ depends only on
         $K = P_z(n)$. Therefore, the expected value of \eqref{eq:preteri}
         is
         \begin{equation}\label{eq:lolkin}O\Big(\mathbb{E}\Big(
         \left(\frac{\log 2 P_z(n)}{\log z} + 1\right)
  \sum_{\substack{d' \mid n:\; \text{$d'$ squarefree} \\ p \mid d' \Rightarrow z^{1/2} < p \leq z}} 1
         \Big)\Big)
         .\end{equation}
         Now, $\log P_z(n) = \sum_{p|n} \log p$. Let $\xi$ denote the random variable given by $$\xi= \sum_{\substack{d' \mid n:\; \text{$d'$ squarefree} \\
             p \mid d' \Rightarrow z^{1/2} < p \leq z}} 1$$ and let $A_p$ denote the event that $p \mid n$. Then \eqref{eq:lolkin}
         is at most a constant times
         \begin{equation}\label{eq:rorkin}
        \mathbb{E}\Big( \xi \Big) +
        \frac{1}{\log z} \sum_{p\leq z} \frac{\log p}{p}
        \mathbb{E}\Big( \xi
         \Big|\; A_p \Big).\end{equation}
         Clearly
         \[\begin{aligned}\mathbb{E}&\Big(\xi
\Big)\leq \frac{1}{N} \sum_{n \leq 2N} \sum_{\substack{d' \mid n:\; \text{$d'$ squarefree} \\ p \mid d' \Rightarrow z^{1/2} < p \leq z}} 1\\
&\ll \frac{1}{N} \sum_{\substack{\text{$d$ square-free} \\ p \mid d \Rightarrow z^{1/2} < p \leq z}} \frac{N}{d} = \sum_{\substack{\text{$d$ square-free} \\ p \mid d \Rightarrow z^{1/2} < p \leq z}} \frac{1}{d} = \prod_{z^{1/2} < p \leq z}\left(1 + \frac{1}{p}\right) \sim \frac{\log z}{\log z^{1/2}} \ll 1.
         \end{aligned}\]
    We must also estimate the conditional expectation:
         for $p\leq z\leq N$,
         \[\begin{aligned}
         \mathbb{E}&\Big(\xi         \Big|\; A_p \Big)
         \ll \frac{1}{N/p} \mathop{\sum_{n \leq 2N}}_{p|n} \sum_{\substack{d' \mid n:\; \text{$d'$ squarefree} \\ p' \mid d' \Rightarrow z^{1/2} < p' \leq z}} 1\\
         &\ll \frac{1}{N/p} \left(\sum_{\substack{\text{$d$ square-free}: p\nmid d \\ p' \mid d \Rightarrow z^{1/2} < p' \leq z}} \frac{N/p}{d} +
         \sum_{\substack{\text{$d$ square-free}: p|d \\ p' \mid d \Rightarrow z^{1/2} < p' \leq z}} \frac{N/p}{d/p}\right)\\ &\ll
        \sum_{\substack{\text{$d$ square-free}: p\nmid d \\ p' \mid d \Rightarrow z^{1/2} < p' \leq z}} \frac{1}{d} \leq 
        \prod_{z^{1/2} < p \leq z}\left(1 + \frac{1}{p}\right) \ll 1.
         \end{aligned}\]

         Hence, the expression in \eqref{eq:rorkin} is
         \[\ll
         1 + \frac{1}{\log z} \sum_{p\leq z} \frac{\log p}{p} \ll
         1 + \frac{\log z}{\log z} \ll 1.\]
\end{proof}

\begin{prop}\label{prop:askesis}
  Define $D(n;a) = \sum_{d|n: d\leq a} \mu(d)$.
  Let $N,A\geq 1$.
  For each $N<n\leq 2 N$, let $A \leq a(n) \leq 2 A$.
  Then, given  the factorization $n= p_1^{\alpha_1} p_2^{\alpha_2} \dotsb p_r^{\alpha_r}$, where $p_1<p_2<\dotsc<p_r$,
  Algorithm \ref{alg:factofun} computes $D(n;a(n))$.
  in expected time $O(\log \log N)$ on average over $n=N+1,\dotsc, 2 N$.
\end{prop}
\begin{proof}
  Algorithm \ref{alg:factofun} computes $D(n;a)$ recursively:
  it calls itself to compute $D(n_0;a)$ and $D(n_0;a/p_r)$, where
  $n_0 = p_1 p_2 \dotsb p_{r-1}$, 
  and then returns $D(n;a) = D(n_0;a) - D(n_0;a/p_r)$.
  The contribution of $D(n_0;a)$ is that of divisors $\ell|n$ with $p_r\nmid \ell$,
  whereas the contribution of $D(n_0;a/p_r)$ corresponds to that
  of divisors $\ell|n$ with $p_r|\ell$.
  
  The algorithm terminates in any of three circumstances:
  \begin{enumerate}
  \item\label{it:dada1} for $a<1$, returning $D(n;a)=0$,
  \item\label{it:dada2} for $n=1$ and $a\geq 1$, returning $D(n;a)=1$,
  \item\label{it:dada3} for $n>1$ and $a\geq n$, returning $D(n;a)=0$.
  \end{enumerate}
  Here it is evident that the algorithm gives the correct output for the cases \eqref{it:dada1}--\eqref{it:dada2}, whereas
  case \eqref{it:dada3} follows from $D(n;a) =
   \sum_{d|n: d\leq a} \mu(d) = \sum_{d|n}\mu(d) = 0$ for $n>1$, $a\geq n$.

   We can see recursion as traversing a {\em recursion tree}, with
   leaves corresponding to cases in which the algorithm terminates. (In the study of algorithms, trees
   are conventionally drawn with the root at the top.)
   The total running time is proportional to the number of vertices in the tree.
   If the algorithm were written to terminate only for $n=1$, the tree
   would have $2^r$ leaves; as it is, the algorithm is written so that
   some branches terminate at depth much lower than $r$. 
   We are to bound the average number of vertices of the recursion tree
   for inputs $N<n\leq 2 N$ and $a=a(n) \in [A,2A]$.

   Say we are at the depth reached after taking care of all $p_i$ with
   $p_i>z$. The branches that have survived correspond to
   $d|n$ with $p|d \Rightarrow p>z$, $d\leq 2 A$ and $d>A/P_z(n)$.
   We are to compute $D(P_z(n);a/d)$. (If $d>2 A$, then $a/d<1$, and
   so our branch has terminated by case \eqref{it:dada1} above.
   If $d\leq A/P_z(n)$, then $a/d\geq P_z(n)$, and we are in case
   \eqref{it:dada3}.) 

   Now we continue running the algorithm until we take care of all $p_i$
   with $p_i>z^{1/2}$. On each branch that survived up to depth $p>z$, the vertices between that depth and depth $p>z^{1/2}$ correspond to square-free divisors
   $d'|n$ such that $p|d\Rightarrow z^{1/2}<p\leq z$.

   By Lemma \ref{lem:isis}, we conclude that the average number of nodes
   in the tree corresponding to $z^{1/2}<p\leq z$ is $O(1)$. Letting
   $z = N, N^{1/2}, N^{1/4}, N^{1/8},\dotsc$, we obtain our result.

\end{proof}



In this way, letting $\Delta = \sqrt{x/v}$,
we can compute $D(r;x/r)$ for all $x/u<r\leq x/v$
in time $O((x/v) \log \log(x/v))$ and space $O(\sqrt{x/v} \log(x/v))$.
Summing values of $D(r;x/r)$ for successive values of $r$ to compute
$S(m) = \sum_{r\leq m} D(r;x/r)$ for $x/u<m\leq x/v$ takes 
time $O(x/v)$ and additional space\footnote{One may take a little more space (but no more than $O(\sqrt{x/v} \log(x/v))$)
  if one decides to parallelize this summation procedure.} $O(1)$.
As $a$ decreases and $m=x/a$ increases, we
may (and should) discard values of $S(m)$ and $D(r;x/r)$ that we no longer
need, so as to keep space usage down.

We have thus shown that we can compute
the right side
of (\ref{eq:pina})
 in time $O((x/v) \log \log x)$ and space $O(\sqrt{x/v}\cdot \log x)$ for any
$1\leq v\leq u = \sqrt{x}$. 

It is easy to see that, if we use the algorithm in
\cite[Main Thm.]{Helfgott2020} instead of the classical segmented sieve of
Eratosthenes, we can accomplish the same task in time
$O((x/v) \log x)$ and space $O((x/v)^{1/3} (\log x)^{5/3})$.

{\bf A few words on the implementation.} See Algorithm \ref{alg:largvar}.

{\em Choice of $\Delta$.}
The size of the segments used by the sieve is to be chosen at the outset:
 $\Delta = C \max(\sqrt{u},\sqrt{x/v}) = C \sqrt{x/v}$
  (for some choice of constant $C\geq 1$) if we use
 the classical segmented sieve (\textsc{SegFactor}), or
 \begin{equation}\label{eq:betdelt}\Delta = C \max\left(\sqrt[3]{u} (\log u)^{2/3}, \sqrt[3]{\frac{x}{v}} (\log x/v)^{2/3}\right)
 = C \sqrt[3]{\frac{x}{v}} \left(\log \frac{x}{v}\right)^{2/3}\end{equation} for
 the improved segmented sieve in \cite[Main Thm.]{Helfgott2020}.


{\em Memory usage.} It is understood that calls such
as $F\gets \textsc{SegFactor}(a_0,\Delta)$ will result in freeing or reusing
the memory previously occupied by $F$. (In other words,
 ``garbage-collection'' will be taken
 care of by either the programmer or the language.)

{\em Parallelization.}  Most of the running time is spent in function
 \textsc{SArr} (Algorithm \ref{alg:sarr}), which is easy to parallelize.
 We can let each processor sieve a block of length $\Delta$.
 Other than that -- the issue
 of computing an array of sums $\mathbf{S}$
 (as in Algorithm \ref{alg:sarr}) in parallel
 is a well-known problem ({\em prefix sums}), for which solutions of
 varying practical efficiency are known. We follow a common two-level
 algorithm: first, we divide the array into as many blocks as there are
 processing elements; then (level 1) we let each processing element compute,
 in parallel,
 an array of prefix sums for each block, ending with the total of the block's
 entries; then we compute prefix sums
 of these totals to create offsets; finally (level 2), we let each
 processing element add its block's offset to all elements of its block.

 
\section{The case of a large free variable}\label{sec:largefree}

We now show how to compute the first double sum on the righthand side of \eqref{eq:nemero}. That double sum equals
\begin{align}\label{eq:largefreesum} \sum_{m, n \leq v} \mu(m) \mu(n) \left\lfloor \frac{x}{mn} \right\rfloor.\end{align} Note that, in \cite{DelegliseRivat96}, this turns out to be the easy case. However, they take $v = x^{1/3}$, while we will take $v = x^{2/5}$. As a result, we have to take much greater care with the computation to ensure that the run time does not become too large. 

\subsection{A first try}
We begin by splitting $[1, v] \times [1, v]$ into neighborhoods $U$ around points $(m_0, n_0)$. For simplicity, we will take these neighborhoods to be rectangles of the form $I_x \times I_y$ with $I_x = [m_0 - a, m_0 + a)$ and $I_y = [n_0 - b, n_0 + b)$, where $\sqrt{m_0}\ll a < m_0$ and $\sqrt{n_0}\ll b < n_0$. (In Section \ref{sec:optimalchoice}, we will partition the two intervals $[1,v]$ into intervals of the form $[x_0, (1+\eta)x_0)$ and $[y_0, (1+\eta)y_0)$, with $0< \eta \leq 1$ a constant. We will then specify $a$ and $b$ for given $x_0$ and $y_0$, and subdivide $[x_0, (1+\eta)x_0) \times [y_0, (1+\eta)y_0)$ into rectangles $I_x \times I_y$ with $|I_x| = 2 a$
            and $|I_y| = 2 b$.) Applying a local linear approximation to the function $\frac{x}{mn}$ on each neighborhood yields
            \begin{equation}\label{eq:beata}\frac{x}{mn} = \frac{x}{m_0 n_0} + c_x(m - m_0) + c_y (n - n_0) + \mathrm{ET}_{\mathrm{quad}}(m,n),\end{equation}
            where $\mathrm{ET}_{\mathrm{quad}}(m,n)$ is a quadratic error term (that is, a term whose size is bounded by $O(\max(n-n_0, m-m_0)^2)$ and $$c_x = \frac{-x}{m_0^2 n_0}, \ c_y = \frac{-x}{m_0 n_0^2}.$$ The quadratic error term will be small provided that $U$ is small. We will show how to choose $U$ optimally at the end of this section. The point of applying the linear approximation is that it will ultimately allow us to separate the variables in our sum. The one complicating factor is the presence of the floor function. If we temporarily ignore both the floor function in \eqref{eq:largefreesum} and the quadratic error term, we can see very clearly how the linear approximation helps us.
            To wit:
            \begin{equation}\label{eq:wouldtake}
              \sum_{(m, n) \in I_x \times I_y} \mu(m) \mu(n) \frac{x}{mn}
            \end{equation}
            is approximately equal to 
\begin{align}&\sum_{(m, n) \in I_x \times I_y} \mu(m)\mu(n) \left(\frac{x}{m_0 n_0} + c_x(m - m_0) + c_y(n - n_0)\right) \notag \\ &= \left(\sum_{m \in I_x} \mu(m) \left(\frac{x}{m_0 n_0} + c_x(m - m_0)\right)\right) \cdot \sum_{n \in I_y} \mu(n) \notag \\ &+ \left(\sum_{n \in I_y} \mu(n) c_y (n - n_0)\right) \cdot \sum_{m \in I_x} \mu(m). \label{eq:nofloorsum} \end{align} 

One can use the segmented sieve of Eratosthenes 
to compute the values of $\mu(m)$ for $m \in I_x$ and $\mu(n)$ for $n \in I_y$. If $a < \sqrt{x_0}$ or $b < \sqrt{y_0}$, we compute the values of $\mu$ in segments of length about $\sqrt{x_0}$ or $\sqrt{y_0}$ and use them for several neighborhoods $I_x \times I_y$. In any event, computing \ref{eq:nofloorsum} given $\mu(m)$ for $m \in I_x$ and $\mu(n)$ for $n \in I_y$ takes only time $O(\max(a, b))$ and negligible space. 

\subsection{Handling the difference between reality and an approximation}\label{subs:handl}

Proceeding as above, we can compute
the sum
\[ S_0:= \sum_{(m, n) \in I_x \times I_y} \mu(m) \mu(n) \left(\left\lfloor \frac{x}{m_0 n_0} + c_x(m - m_0)\right\rfloor + \left\lfloor c_y(n - n_0)\right\rfloor \right)\]
in time $O(\max(a,b))$ and space $O(\log \max(x_0,y_0))$,
given arrays with the values of $\mu(m)$ and $\mu(n)$.
The issue is that $S_0$ is not the same as
\begin{align*} S_1 := \sum_{(m, n) \in I_x \times I_y} \mu(m) \mu(n) \left(\left\lfloor \frac{x}{m_0 n_0} + c_x(m - m_0) + c_y(n - n_0)\right\rfloor
  \right),\end{align*} and it is certainly not the same as the sum we actually want to
compute, namely,
\begin{align*} S_2 := \sum_{(m, n) \in I_x \times I_y} \mu(m) \mu(n)
  \left\lfloor \frac{x}{m n}\right\rfloor .
\end{align*}

From now on, we will write
\[L_0(m,n) = \left\lfloor\frac{x}{m_0 n_0} + c_x(m - m_0)\right\rfloor + \left\lfloor c_y(n - n_0)\right\rfloor,\]
\[L_1(m,n) = \left\lfloor \frac{x}{m_0 n_0} + c_x(m - m_0) + c_y(n - n_0)\right\rfloor,\;\;\;\;\; L_2(m,n) = \left\lfloor \frac{x}{m n}\right\rfloor.\]
Here $m_0$, $n_0$ and $x$ are understood to be fixed.
Our challenge will be to show that the weights $L_2-L_1$ and $L_1-L_0$
actually have a simple form -- simple enough that $S_2-S_1$ and $S_1-S_0$
can be computed quickly.

We approximate $c_y$ by a rational number $a_0/q$ with $q \leq Q = 2b$
such that $$\delta:= c_y - a_0/q$$ satisfies
$|\delta| \leq 1/q Q.$
Thus, \begin{equation}\label{eq:nueve}
  \left|c_y (n-n_0) - \frac{a_0 (n-n_0)}{q}\right|\leq \frac{1}{2 q}
.\end{equation}
We can find
such an $\frac{a_0}{q}$ in time $O(\log Q)$ using continued fractions
(see Algorithm \ref{alg:diophfrac}).

Write $r_0 = r_0(m)$ for the integer such that the absolute value of
\begin{equation}\label{eq:ultima}
  \beta = \beta_m :=
  \left\{\frac{x}{m_0 n_0} + c_x (m-m_0)\right\} - \frac{r_0}{q}
        \end{equation}
          is minimal (and hence $\leq 1/2 q$). If there are two such
          values, choose the greater one. Then
          \begin{equation}\label{eq:reinas}
            - \frac{1}{2 q} \leq \beta< \frac{1}{2 q}
            .\end{equation}

          We will later make sure that we choose our neighborhoods
          $I_x\times I_y$ so that
          $|\mathrm{ET}_{\mathrm{quad}}(m,n)|\leq 1/2 b$, where
          $\mathrm{ET}_{\mathrm{quad}}(m,n)$ is defined by \eqref{eq:beata}.
          We also know that $\mathrm{ET}_{\mathrm{quad}}(m,n)>0$, since the
          function
          $(m,n)\mapsto x/m n$ is convex. 
          We are of course assuming that $I_x\times I_y$ is contained in the
          first quadrant, and so           $(m,n)\mapsto x/m n$ is
          well-defined on it.

          The aforementioned notation will be used throughout this section.
          
          \begin{lem}\label{lem:tofutti}
            Let $(m,n)\in I_x\times I_y$.
            Unless $a_0(n - n_0) + r_0 \in \{0, -1\} \mo q$, 
            \begin{equation*}
              L_2(m,n) = L_1(m,n).
            \end{equation*}
          \end{lem} 
          
\begin{proof}
Since $0 < \mathrm{ET}_{\mathrm{quad}}(m, n) \leq 1/2b$, we can have
          \begin{equation}\label{eq:sulliv}
          \left\lfloor \frac{x}{m n}\right\rfloor \ne
          \left\lfloor \frac{x}{m_0 n_0} + c_x(m - m_0) + c_y(n - n_0)\right\rfloor\end{equation}
          (in which case the left side equals the right side plus $1$)
          only if
          \begin{equation}\label{eq:ierro}\left\{
          \frac{x}{m_0 n_0} + c_x(m - m_0) + c_y(n - n_0)\right\} \geq
          1 - \frac{1}{2 b}.\end{equation} 
          Since $q\leq 2 b$ and
          \[
          \frac{x}{m_0 n_0} + c_x(m - m_0) + c_y(n - n_0)
\in   \frac{a_0 (n-n_0) + r_0}{q} + \left[-\frac{1}{q},\frac{1}{q}\right),\]
          we see that \eqref{eq:ierro} can be the case only if
          $a_0 (n-n_0) + r_0$ is in $\{0,-1\} \mo q$.

\end{proof}         
         
\begin{lem}\label{lem:cuy}
              Let $(m,n)\in I_x\times I_y$.
Unless $a_0(n - n_0) + r_0 \equiv 0 \pmod{q}$, 
\begin{align}\label{eq:katnip}
L_1(m,n) - L_0(m,n) &=
\begin{cases}
  0 & \text{if $r_0 + \overline{a_0 (n-n_0)} \leq q$},\\
    1 & \text{otherwise,}
\end{cases}\\ &+
\begin{cases}
  1 & \text{if $q|(n-n_0) \wedge (\delta (n-n_0)< 0)$},\\
    0 & \text{otherwise.}
\end{cases} 
\end{align}
\end{lem}  
  
 \begin{proof}
 Recall that, for all real numbers $A$ and $B$, 
\[\lfloor A + B \rfloor - (\lfloor A \rfloor + \lfloor B \rfloor) = 
\begin{cases}
0, & \mathrm{if} \ \{A\} + \{B\} < 1 \\
1, & \mathrm{otherwise.}
\end{cases}
\]
Thus, $L_1(m,n)-L_0(m,n)$
is either $0$ or $1$, and it is $1$ if and only if
\begin{equation}\label{eq:groningen}
\left\{ \frac{x}{m_0 n_0} + c_x(m - m_0)\right\} + \left\{ c_y(n - n_0)\right\}
\end{equation}
is $\geq 1$. By \eqref{eq:nueve} and \eqref{eq:reinas},
the quantity in \eqref{eq:groningen}
lies in
\[\frac{r_0}{q} + \left\{\frac{a_0 (n-n_0)}{q}\right\} 
+ \left[-\frac{1}{q}, \frac{1}{q}\right)\]
  unless, possibly, if $a_0 (n-n_0)\equiv 0 \mo q$, that is, if
  $q|(n-n_0)$.
 Hence, unless $a_0 (n- n_0) + r_0 \equiv 0 \mo q$ or $q|(n-n_0)$,
  the expression in \eqref{eq:groningen} is $\geq 1$ if and only if
  $r_0/q + \{a_0 (n-n_0)/q\} \geq 1$.
  Moreover, if $q|(n-n_0)$ but
  $a_0 (n- n_0) + r_0 \not\equiv 0 \mo q$,
  it is easy to see that the expression in \eqref{eq:groningen} is $<1$
  iff $\delta (n-n_0) = c_y (n-n_0) - a_0 (n-n_0)/q$ is $\geq 0$.
 \end{proof} 

 It follows immediately from Lemmas \ref{lem:tofutti} and \ref{lem:cuy}
 that
  \begin{equation}\label{eq:adaman}
    L_2(m,n)-L_0(m,n)  = \begin{cases}
    0 & \text{if $r_0 + \overline{a_0 (n-n_0)} \leq q$,}\\
    1 & \text{otherwise,}
    \end{cases}
  \end{equation}
  unless $r_0 + a_0 (n-n_0) \in \{0,-1\} \mo q$,
  where we write $\overline{a}$ for the integer in $\{0,1,\dotsc,q-1\}$
  congruent to $a$ modulo $q$.
  
 Note that the first term on the right side of
 \eqref{eq:adaman} depends only on $n \mo q$ (and $a_0 \mo q$ and $r_0$),
 and the second term depends only on $n\mo q$, $\sgn(n-n_0)$ and $\sgn(\delta)$
 (and not on $r_0$; hence it is independent of $m$). 
Given the values of $\mu(n)$ for $n\in I_y$,
  it is easy to make a table of
  \[\rho_r = \mathop{\sum_{n\in I_y}}_{a_0 (n-n_0)\equiv r \mo q} \mu(n)\]
  for $r\in \mathbb{Z}/q \mathbb{Z}$ in time $O(b)$
  and space $O(q \log b)$,
  and then a table of
  \[\sigma_r = \mathop{\sum_{n\in I_y}}_{\overline{a_0 (n-n_0)} > q- r} \mu(n)\]
  for $0\leq r\leq q$ in time $O(q)$ and space $O(q \log b)$.
  We also compute
  \[\mathop{\mathop{\sum_{n\in I_y}}_{q|n-n_0}}_{\delta \cdot (n-n_0) < 0} \mu(n)\]
  once and for all.   It remains to deal with the problematic cases 
  $a_0 (n-n_0)+r_0\in \{0,-1\} \mo q$.

    \begin{lem}\label{lem:tannutuva}
      Let $(m,n)\in I_x\times I_y$.
 If $a_0(n - n_0) + r_0 \equiv -1 \pmod{q}$ and $q>1$, then 
 \begin{equation*}
L_2(m,n) - L_1(m,n)
   = \begin{cases}
    1 & \text{if $n \not\in I$,}\\
    0 & \text{if $n \in I$,}
    \end{cases}
    \end{equation*}
 where $I = (\mathbf{x}_-,\mathbf{x}_+)$ if
  the equation \[\gamma_2 \mathbf{x}^2 + \gamma_1 \mathbf{x} + \gamma_0 = 0\]
  has real roots $\mathbf{x}_-< \mathbf{x}_+$, and $I=\emptyset$
  otherwise. Here $\gamma_0 = x q$, $\gamma_2 = - a_0 m$
  and
   \begin{align*}
    \gamma_1 &= 
    \left( -\left\lfloor \frac{x}{m_0 n_0} + c_x(m - m_0)\right\rfloor q -
    (r_0+1) + a_0 n_0\right) m .
  \end{align*}

 \end{lem} 
  
  \begin{proof} 
The question is whether $L_2(m,n)>L_1(m,n)$.
Since \begin{equation}\label{eq:tofucat}  -1/2q \leq \beta<1/2 q \ \mathrm{and} \ 
|\delta (n-n_0)|\leq 1/2 q,\end{equation}
we know that
\begin{align*}
  &\left\{\frac{x}{m_0 n_0} + c_x(m - m_0) + c_y(n - n_0)\right\}
= \left\{\frac{r_0}{q} + \beta + \frac{a_0 (n-n_0)}{q} + \delta (n-n_0)\right\}
  \\ &= \left\{- \frac{1}{q} + \beta + \delta (n-n_0)\right\}
  = \frac{q-1}{q} + \beta + \delta (n-n_0),
\end{align*}
where the last line follows from \eqref{eq:tofucat}. 
Hence, $L_2(m,n)>L_1(m,n)$ if and only if
\begin{equation}\label{eq:wookie}\frac{x}{m n} - \left(
\frac{x}{m_0 n_0} + c_x(m - m_0)
+ c_y(n - n_0)\right)
  \geq \frac{1}{q} - \beta - \delta (n-n_0).\end{equation}
  This, in turn, is equivalent to
  \begin{equation}\label{eq:gunth}
    \frac{c_0}{n} + c_1 + c_2 n \geq 0,\end{equation}
  where $c_0 = x/m$, $c_2 = -a_0/q$ and
  \begin{align*}
    c_1 &= -\left(\frac{x}{m_0 n_0} + c_x(m - m_0)-\beta\right)
  + \frac{a_0}{q} n_0 - \frac{1}{q} \\ &=
  -\left\lfloor \frac{x}{m_0 n_0} + c_x(m - m_0)\right\rfloor -
  \frac{r_0+1}{q} + \frac{a_0}{q} n_0 
  .\end{align*}

  Since $a_0/q$ is a Diophantine approximation to $c_y = - x/m_0 n_0^2 < 0$,
  it is clear that $a_0/q$ is non-positive. Consequently, if $q>1$,
  $a_0$ must be negative, since $a_0$ and $q$ are coprime.
  Hence, $c_2$ is positive, and so \eqref{eq:gunth} holds iff
  $n\not\in I$, where $I = (\mathbf{x}_-,\mathbf{x}_+)$ if
  the equation \[c_2 \mathbf{x}^2 + c_1 \mathbf{x} + c_0 = 0\]
  has real roots $\mathbf{x}_-\leq\mathbf{x}_+$, and $I=\emptyset$
  otherwise. 
  
  
  \end{proof}
  
  Solving a quadratic equation is not computationally expensive;
  in practice, the function $x\mapsto \lfloor \sqrt{x}\rfloor$ generally
  takes less time to compute than a division. Thus it makes sense to
  consider it to take $O(1)$ time, since we are thinking of the four basic
  operations as taking $O(1)$ time.

  What we have to do is keep a table of
  \[\rho_{r,\leq n'} = \mathop{\sum_{n\in I_y, n\leq n'}}_{a_0 (n-n_0)\equiv r \mo q}
  \mu(n).\]
  We need only consider values of $n'$ satisfying $a_0 (n' - n_0)\equiv r \mo q$
  (since $\rho_{r,\leq n'} = \rho_{r,\leq n''}$ for $n''$ the largest number
  $n''\leq n'$ with $a_0 (n'' - n_0)\equiv r \mo q$).
  It is then easy to see that we can construct the table in time
  $O(b)$ and space $O(b \log b)$, simply letting $n$ traverse $I_y$ from left
  to right. (In the end, we obtain $\rho_r$ for every $r\in \mathbb{Z}/
  q\mathbb{Z}$.) In the remaining lemmas, we show how to handle the cases where $a_0(n - n_0) + r_0 \equiv 0 \pmod{q}$.

  \begin{lem}\label{lem:guineapigs!}       Let $(m,n)\in I_x\times I_y$.
    If $a_0(n - n_0) + r_0 \equiv 0 \pmod{q}$, then 
    \begin{align*}
L_1(m,n) - L_0(m,n) = 
\begin{cases}
    0 & \text{if $n \not\in I$,}\\
    1 & \text{if $n \in I$,}
    \end{cases}
    \end{align*}
    where, if $r_0\not\equiv 0 \mo q$,
    \[I = \begin{cases}
n_0 - \frac{\beta}{\delta} + \frac{1}{\delta} \cdot [0,\infty)
      &\text{if $\delta \ne 0$,}\\
    \mathbb{R} & \text{if $\delta=0$ and $\beta\geq 0$,}\\
    \emptyset & \text{if $\delta=0$ and $\beta < 0$,}\end{cases}\]
    and, if $r_0\equiv 0 \mo q$,
    \[I = \begin{cases} \mathbb{R} &\text{if $\beta<0$ and $\delta<0$}\\
(-\infty,n_0] \cup [n_0 - \frac{\beta}{\delta},\infty)
      &\text{if $\beta <0$ and $\delta>0$}\\
        n_0 + \frac{1}{\delta} [-\beta,0)
        &\text{if $\beta>0$ and $\delta\ne 0$,}\\
    \emptyset & \text{otherwise.}\end{cases}\]
    \end{lem}
\begin{proof}
Since $\{a_0 (n-n_0)/q\} = \{-r_0/q\}$,
\begin{align*}
  \left\{ \frac{x}{m_0 n_0} + c_x(m - m_0)\right\} + \left\{ c_y(n - n_0)\right\}
  &= \left\{ \frac{r_0}{q} + \beta\right\} +
 \left\{ -\frac{r_0}{q} + \delta (n- n_0)\right\}.
\end{align*}
Recall that  $-1/2q\leq \beta<1/2 q$ and $|\delta (n-n_0)|\leq 1/2 q$.
For $r_0\not\equiv 0 \mo q$,
$\{r_0/q+\beta\}+\{-r_0/q + \delta (n-n_0)\}\geq 1$ iff
$\beta+\delta(n-n_0)\geq 0$.
We treat the case $r_0 \equiv 0 \mo q$ separately:
$\{\beta\}+\{\delta (n-n_0)\}\geq 1$
iff either (a) $\beta<0$ and $\delta (n-n_0)<0$, or
(b) $\beta \delta (n-n_0) < 0$ and $\beta+\delta(n-n_0)\geq 0$.

\end{proof}

\begin{lem}\label{lem:felix}    Let $(m,n)\in I_x\times I_y$.
  If $a_0(n - n_0) + r_0 \equiv 0 \pmod{q}$ and $q>1$,
  \begin{equation*}     L_2(m,n) - L_1(m,n) = 
 \begin{cases}
    0 & \text{if $n \not\in I\cap J$,}\\
    1 & \text{if $n \in I\cap J$,}
    \end{cases}
    \end{equation*}
where $I = [\mathbf{x}_-,\mathbf{x}_+]$ if
  the equation \[\gamma_2 \mathbf{x}^2 + \gamma_1 \mathbf{x} + \gamma_0 = 0\]
  has real roots $\mathbf{x}_-\leq \mathbf{x}_+$, and $I=\emptyset$
  otherwise, whereas
  $J = n_0 -\beta/\delta -\frac{1}{\delta} (0,\infty)$ if $\delta\ne 0$,
  $J= \emptyset$ if $\delta=0$ and $\beta\geq 0$ and
  $J= (-\infty,\infty)$ if $\delta=0$ and $\beta<0$.
Here $\gamma_0 = x q$, $\gamma_2 = - a_0 m$
  and
   \begin{align*}
    \gamma_1 &= 
    \left( -\left\lfloor \frac{x}{m_0 n_0} + c_x(m - m_0)\right\rfloor q -
    r_0 + a_0 n_0\right) m.
   \end{align*}
\end{lem}

\begin{proof}

As in the proof of Lemma \ref{lem:tannutuva}, we have
\begin{align*}
  \left\{\frac{x}{m_0 n_0} + c_x(m - m_0) + c_y(n - n_0)\right\}
&= \left\{\frac{r_0}{q} + \beta + \frac{a_0 (n-n_0)}{q} + \delta (n-n_0)\right\}
  \\ &= \left\{\beta + \delta (n-n_0)\right\},
\end{align*}
where the last equality follows from the fact that $a_0(n-n_0) + r_0 \equiv 0 \pmod{q}.$
We know that $\beta + \delta (n-n_0)<1/q$, whereas
$0<\mathrm{ET}_{\mathrm{quad}}(m,n)\leq 1/2 b \leq 1/q$.
Since $q>1$, we see that, if  $\beta + \delta (n-n_0) \geq 0$,
the inequality
\begin{equation}\label{eq:dindin}
    \left\lfloor \frac{x}{m n}\right\rfloor >
          \left\lfloor \frac{x}{m_0 n_0} + c_x(m - m_0) + c_y(n - n_0)\right\rfloor\end{equation}
cannot hold.
If $\beta + \delta (n-n_0) < 0$,
then \eqref{eq:dindin} holds iff
\begin{equation}\label{eq:armado}\frac{x}{m n} - \left(\frac{x}{m_0 n_0} + c_x(m - m_0) + c_y(n - n_0)\right)
\geq - \beta - \delta (n-n_0),
\end{equation}
Much as in the proof of Lemma \ref{lem:tannutuva}, this inequality holds iff $n\in I$, where $I=
[\mathbf{x}_-,\mathbf{x}_+]$ if the equation
$c_2 \mathbf{x}^2 + c_1 \mathbf{x} + c_0 = 0$ has real roots
$\mathbf{x}_-\leq \mathbf{x}_+$, where $c_0=x/m$, $c_2 = -a_0/q$ and
\begin{align*}
    c_1 =
  -\left\lfloor \frac{x}{m_0 n_0} + c_x(m - m_0)\right\rfloor -
  \frac{r_0}{q} + \frac{a_0}{q} n_0 
  ,\end{align*}
and $I=\emptyset$ if the equation has complex roots.
\end{proof}

\begin{lem}\label{lem:wuwu}    Let $(m,n)\in I_x\times I_y$.
  If $q=1$,
  \[ L_2(m,n) - L_1(m,n) =     \begin{cases} 0 & \text{if $n \not\in
    (I_0\cap J) \cup (I_1\cap (\mathbb{R}\setminus J))$,}\\
    1 & \text{if $n \in (I_0\cap J) \cup (I_1\cap (\mathbb{R}\setminus J))$,}
  \end{cases}
  \]
  where $J = n_0 -\beta/\delta -\frac{1}{\delta} (0,\infty)$ if $\delta\ne 0$,
  $J= \emptyset$ if $\delta=0$.

  If $a\ne 0$, then
  $I_j = [\mathbf{x}_{-,j},\mathbf{x}_{+,j}]$ if
  the equation \[\gamma_2 \mathbf{x}^2 + \gamma_{1,j} \mathbf{x} + \gamma_0 = 0\]
  has real roots $\mathbf{x}_{-,j}\leq \mathbf{x}_{+,j}$, and $I=\emptyset$
  otherwise.
  Here $\gamma_0 = x q$, $\gamma_2 = - a_0 m$
  and
   \begin{align*}
    \gamma_{1,j} &= 
    \left( -\left\lfloor \frac{x}{m_0 n_0} + c_x(m - m_0)\right\rfloor q -
    (r_0+j) + a_0 n_0\right) m.
   \end{align*}
If $a=0$, then
      \[I_j = \left(-\infty,
  \frac{x}{m} \left(
 \left\lfloor \frac{x}{m_0 n_0} + c_x(m - m_0)\right\rfloor 
    + r_0 + j
  \right)^{-1}
  \right].\]
  \end{lem}
  \begin{proof}
    Just as in the proof of Lemma \ref{lem:felix},
    \[  \left\{\frac{x}{m_0 n_0} + c_x(m - m_0) + c_y(n - n_0)\right\}
    = \left\{\beta + \delta (n-n_0)\right\}.\]
    If $\beta + \delta (n-n_0)<0$, then $L_2(m,n)-L_1(m,n)>0$ holds iff
    \eqref{eq:armado} holds. The term $\delta (n-n_0)$
    cancels out, and so, by \eqref{eq:ultima},
    we obtain that \eqref{eq:armado} holds iff
    \[\frac{x}{m n} \geq
    \left\lfloor \frac{x}{m_0 n_0} + c_x(m - m_0)\right\rfloor 
    + a_0 (n-n_0)
    + r_0,\]
    just as in Lemma \ref{lem:felix}.
    If $\beta + \delta (n-n_0)\geq 0$,  $L_2(m,n)-L_1(m,n)>0$ holds iff
    \eqref{eq:wookie} holds. Again, the term involving $\delta (n-n_0)$
    cancels
    out fully, and so \eqref{eq:armado} holds iff
        \[\frac{x}{m n} \geq
        \left\lfloor \frac{x}{m_0 n_0} + c_x(m - m_0)\right\rfloor 
        + a_0 (n-n_0)
    + r_0 + 1.\]
    \end{proof}

In summary: for a neighborhood $I_x\times I_y$ small enough that
$|\mathrm{ET}_{\mathrm{quad}}(m,n)|\leq 1/2 b$, we need to prepare tables
(in time $O(b)$ and space $O(b \log b)$), compute
a Diophantine approximation (in time $O(\log b)$), and then, for each value of
$m$, we need to (i) compute $r_0=r_0(m)$, (ii) look up
$\sigma_{r_0}$ in a table, (iii) solve a quadratic equation
to account for the case $a_0 (n-n_0)+r_0\equiv -1 \mo q$, 
(iv) solve a quadratic equation and also a linear equation to account
for the case $a_0 (n-n_0) +r_0 \equiv 0 \mo q$.
If $q=1$, then (iii) and (iv) are replaced by the simple task of
computing the expressions in Lemma \ref{lem:wuwu}.
In any event, these are a bounded number of operations
taking a bounded amount of time. Thus, the computation over the neighborhood $I_x \times I_y$ takes total time
$O(a+b)$ and space $O(b \log b)$, given the values of $\mu(m)$ and $\mu(n)$.

\section{Parameter choice. Final estimates.}\label{sec:optimalchoice} What remains now is to choose our neighborhoods $U = I_x \times I_y$ optimally (within a constant factor),
and to specify our choice of $v$. Recall that $I_x = [m_0-a,m_0+a)$,
  $I_y = [n_0-b,n_0+b)$.

    \subsection{Bounding the quadratic error term. Choosing
      \texorpdfstring{$a$}{a} and \texorpdfstring{$b$}{b}.}\label{subs:estquad}
We can use the formula for the error term bound in a Taylor expansion to obtain an upper bound on the error term. Since
$f:(x,y)\mapsto X/x y$ is twice continuously differentiable for $x,y>0$,
we know that, for $(x,y)$ in any convex neighborhood $U$ of any $(x_0,y_0)$ with
$x_0,y_0>0$,
\[
\frac{X}{x y} = \frac{X}{x_0 y_0} + \frac{\partial f(x_0,y_0)}{\partial x}
(x - x_0) + \frac{\partial f(x_0,y_0)}{\partial y} (y-y_0)
+ \mathrm{ET}_{\mathrm{quad}}(x,y),\]
where the {\em Lagrange remainder term}
$\mathrm{ET}_{\mathrm{quad}}(x,y)$ is given by
\[\begin{aligned}
\mathrm{ET}_{\mathrm{quad}}(x,y) &=
 \frac{1}{2} \frac{\partial^2 f(\xi,\upsilon)}{\partial^2 x} (x-x_0)^2 +
\frac{1}{2} \frac{\partial^2 f(\xi,\upsilon)}{\partial^2 y} (y-y_0)^2 \\ &+
\frac{\partial^2 f(\xi,\upsilon)}{\partial x \partial y} (x-x_0) (y-y_0),
\end{aligned}\]
for some $(\xi,\upsilon)=(\xi(x,y),\upsilon(x,y))\in U$ depending on $(x,y)$.
Working with our neighborhood $U=I_x\times I_y$ of $(x_0,y_0)=(m_0,n_0)$,
we obtain that, for $m\in I_x$ and $n\in I_y$,
$|\mathrm{ET}_{\mathrm{quad}}(m,n)|$ is at most
\begin{align}\label{eq:ET_quad} 
  &\leq
  \frac{X}{m'^3 n'}(m - m_0)^2 + \frac{X}{m'^2 n'^2} (m-m_0) (n-n_0) +
  \frac{X}{m'n'^3}(n - n_0)^2,\end{align} where $m' = \min_{(m, n) \in U} m$ and $n' = \min_{(m, n) \in U} n.$ Hence,
by Cauchy-Schwarz, 
\[
|\mathrm{ET}_{\mathrm{quad}}(m,n)|\leq
\frac{3}{2} \left(\frac{X}{m'^3 n'}(m - m_0)^2 + 
\frac{X}{m'n'^3}(n - n_0)^2\right).\]
(From now on, we will write $x$, as we are used to,
instead of $X$, since there
is no longer any risk of confusion with the variable $x$.)

Recall that we need to choose $I_x$ and $I_y$ so that
$\left|\textrm{ET}_{\textrm{quad}}\right|\leq 1/2 b$.
Since $(m - m_0)^2 \leq a^2$ and $(n - n_0)^2 \leq b^2$,
it is enough to require
that $$\frac{x}{m'^3 n'} a^2 \leq \frac{1}{6 b},\;\;\; \ \frac{x}{m' n'^3} b^2 \leq \frac{1}{6 b}.$$
In turn, these conditions hold for
\[a = \sqrt[3]{\frac{(m')^4}{6 x}},\;\;\;
b = \sqrt[3]{\frac{m' (n')^3}{6 x}}.\]
More generally, if we are given that $m'\geq A$, $n'\geq B$ for some
$A$, $B$, we see that we can set
\begin{equation}\label{eq:trupar}
  a = \sqrt[3]{\frac{A^4}{6 x}},\;\;\;
b = \sqrt[3]{\frac{A B^3}{6 x}}.\end{equation}

At the end of Section \ref{sec:largefree}, we showed that it takes time
$O(a+b)$ and space $O(b \log b)$ for our algorithm to run over each neighborhood $I_x \times I_y$. 
Recall that we are dividing $[1, v] \times [1,v]$ into dyadic boxes (or, at any rate, boxes of the form
$\mathbf{B}(A,B,\eta) = [A,(1+\eta) A) \times [B,(1+\eta) B)$, where $0<\eta\leq 1$
    is a constant)
    and that these boxes are divided into neighborhoods $I_x \times I_y$.
    We have $\ll \frac{A B}{ab}$ neighborhoods $I_x \times I_y$ in the
box $\mathbf{B}(A,B,\eta)$.
Thus, assuming that $A\geq B$,
it takes time $$O\left(\frac{A B}{a b}(a + b)\right) =
O\left(\frac{A B}{b}\right) =
O\left(A^{2/3} x^{1/3}\right)$$ to run over this box, using the values
of $a$ and $b$ in \eqref{eq:trupar}.

Now, we will need to sum over all boxes $\mathbf{B}(A,B,\eta)$.
Each $A$ is of the form $\lceil (1+\eta)^i\rceil$ and each $B$ is of the form
$\lceil (1+\eta)^j\rceil$ for $1 \leq (1+\eta)^i, (1+\eta)^j \leq v.$ By symmetry, we may take 
$j\leq i$, that is, $A\geq B$.
Summing over all boxes takes time
\[\begin{aligned} \ll
\sum_{i:(1+\eta)^i \leq v} \sum_{j\leq i} ((1+\eta)^i)^{2/3} x^{1/3} &\ll
\sum_{i:(1+\eta)^i \leq v} i ((1+\eta)^i)^{2/3} x^{1/3} \\ \ll
(\log v) v^{2/3} x^{1/3}&\leq v^{2/3} x^{1/3} \log x.\end{aligned}\]

We tacitly assumed that $a\geq 1$, $b\geq 1$, and so we need to handle
the case of $a<1$ or $b<1$ separately, by brute force. It actually
makes sense to treat the broader case of $a<C$ or $b<C$ by brute force,
where $C$ is a constant of our choice.
The cost of brute-force summation for $(m,n)$ with
$n\leq m\ll (C^3 x)^{1/4}$ (as is the case when $a<C$)
is
\[\ll ((6 C^3 x)^{1/4})^2 \ll x^{1/2},
\]
whereas the cost of brute-force summation for
$(m,n)$ with $m\leq v$, $n\ll (6 x/m)^{1/3}$
(as is the case when $b<C$) is
\[\ll \sum_{m\leq v} \frac{x^{1/3}}{m^{1/3}} \ll x^{1/3} v^{2/3}.\]

Lastly, we need to take into account the fact that we had to pre-compute a list of values of $\mu$ using a segmented sieve (Algorithm \ref{alg:segsievemu}), which takes time $O(v^{3/2} \log \log x)$ and space $O(\sqrt{v} \log \log v)$. Putting everything together, we see that the large free variable case (Section \ref{sec:largefree}) takes time $O(v^{2/3} x^{1/3} \log x + v^{3/2} \log \log x)$ and space $O(\sqrt{v} \log \log x + (v^4/x)^{1/3} \log x),$ where the space bound comes from substituting $b = \sqrt[3]{\frac{m' (n')^3}{6 x}}$ into the space estimate that we had for each neighborhood and adding it to the space bound from the segmented sieve.

\subsection{Choice of \texorpdfstring{$v$}{v}. Total time and space estimates.}

Recall that the case of a large non-free variable (Algorithm \ref{alg:largvar}) takes time $O((\frac{x}{v} + u) \log \log x)$ and space $O(\sqrt{\max(x/v, u)} \log x)$. At the end of Section \ref{sec:largenonfree}, we took $u=\sqrt{x}$,
making the running time $O(\frac{x}{v} \log \log x)$ and space
$O(\sqrt{x/v} \log x)$.

On the other hand, as we just showed,
the case of a large free variable (Algorithm \ref{alg:smallnonfree}) takes time $O(v^{2/3} x^{1/3} \log x + v^{3/2} \log \log x)$ and space $O(\sqrt{v} \log \log x + (v^4/x)^{1/3} \log x)$.

Thus, in order to minimize our running time, we set the two time bounds equal to one another and solve for $v$, yielding
$v = x^{2/5} (\log \log x)^{3/5}/(\log x)^{3/5}$.
Using this value of $v$ (or
any value of $v$ within a constant factor $c$ of it)
allows us to obtain
  \[\mathrm{time} \ \ O\left(x^{\frac{3}{5}} (\log x)^{\frac{3}{5}}
  (\log \log x)^{\frac{2}{5}} \right) \ \
  \mathrm{and \ space} \ \ O\left(x^{\frac{3}{10}} (\log x)^{\frac{13}{10}}
  (\log \log x)^{-\frac{3}{10}} \right),\]
  as desired. Note that our algorithm for the case of
a large non-free variable uses more memory, by far, than that for the case of
a large free variable.

The constant $c$ can be fine-tuned by the user or programmer. It is actually best to set it so that the time taken by the case of a large free variable and by the case of a large non-free variable are within a constant factor
  of each other without being approximately equal.

  If we were to use \cite{Helfgott2020} to factor integers in
  \textsc{SArr} (Algorithm \ref{alg:sarr}) then \textsc{LargeNonFree}
  (Algorithm \ref{alg:largvar}) would take time
  $O((x/v) \log x)$
  and space $O((x/v)^{1/3} (\log (x/v))^{5/3})$. It would then be best to set
  $v = c\cdot x^{2/5}$ for some $c$, leading to total
  time $O(x^{3/5} \log x)$ and total space
  $O\left(x^{1/5} (\log x)^{5/3}\right)$.
  
\section{Implementation details}\label{sec:impdet}

We wrote our program in C++ (though mainly simply in C).
We used gmp (the GNU MP multiple precision library) for a few operations,
but relied mainly on 64-bit and 128-bit
arithmetic. Some key procedures were parallelized by means of OpenMP pragmas.

{\em Basics on better sieving.}
Let us first go over two well-known optimization techniques. The first one is
useful
for sieving in general; the second one is specific to the use of sieves
to compute $\mu(n)$.

\begin{enumerate}
\item When we sieve (function \textsc{SegPrimes}, \textsc{SegMu} or
  \textsc{SegFactor}), it is useful to first compute how our sieve
  affects a segment of length $M=2^3\cdot 3^2\cdot 5\cdot 7\cdot 11$, say.
  (For instance, if we are sieving for primes, we compute which elements
  of $\mathbb{Z}/M\mathbb{Z}$ lie in $(\mathbb{Z}/M \mathbb{Z})^*$.)
  We can then copy that segment onto our longer segment repeatedly,
  and then start sieving by primes and prime powers not dividing $M$.
\item\label{it:ardar} As is explained in \cite{Kuznetsov11} and
  \cite{Hurst18}, and for that matter in \cite[\S 4.5.1]{Helfbook}:
  in function \textsc{SegMu}, for $n\leq x_0 = n_0 + \Delta$,
  we do not actually need to store $\Pi_j = \sum_{p\leq \sqrt{x_0}: p|n} p$;
  it is enough to store
  $S_j \sum_{p\leq \sqrt{x_0}} \lceil \log_4 p\rceil$. The reason is that
  (as can be easily checked)
  $\Pi_j < \prod_{p|n} p$ if and only if $S_j < \lceil \log_4 n\rceil$.
  In this way, we use space $O(\Delta \log \log x_0)$ instead of space
  $O(\Delta \log x_0)$. We also replace many multiplications by additions;
  in exchange,
  we need to compute $\lceil \log_4 p\rceil$ and $\lceil \log_4 n\rceil$,
  but that takes very little time, as it only involves
  counting the space occupied by $p$ or $n$ in base $2$, and that is a task
  that a processor can usually accomplish extremely quickly. 
  \end{enumerate}
Technique (\ref{it:ardar}) here is not essential in our context,
as \textsc{SegMu} is not a bottleneck, whether for time or for space.
It is more important to optimize factorization -- as we are about to explain.

{\em Factorizing via a sieve in little space.} We wish to store the list of
prime factors of a positive number $n$ in at most twice as much
space as it takes to store $n$. We can do so simply and rapidly as follows.
We initialize $a_n$ and $b_n$ to $0$. When we find a new prime factor $p$,
we reset $a_n$ to $2^k a_n + 2^{k-1}$, where $k = \lfloor \log_2 p\rfloor$,
and $b_n$ to $2^k b_n + p - 2^k$. In the end, we obtain, for example,
\[a_{2\cdot 3\cdot 5\cdot 7} = 111010_2,\;\;
b_{2\cdot 3\cdot 5\cdot 7} = 010111_2.\]
We can easily read the list of prime factors $2$, $3$, $5$, $7$ of
$n=2\cdot 3\cdot 5\cdot 7$
from $a_n$ and $b_n$, whether in ascending or in descending order:
we can see $a_n$ as marking where each prime in $b_n$ begins, as well
as providing the leading $1$: $2=\mathbf{1}0_2$, $3=\mathbf{1}1_2$,
$5=\mathbf{1}01_2$, $7=\mathbf{1}11_2$.

The resulting savings in space lead to a significant speed-up in practice,
due no doubt in part to better cache usage. The bitwise operations required
to decode the factorization of $n$ are very fast, particularly if one is
willing to go beyond the $C$ standard; we used instructions
available in gcc (\texttt{\_\_builtin\_clzl},
\texttt{\_\_builtin\_ctzl}).

{\em Implementing the algorithm in integer arithmetic.}
Manipulating rationals is time consuming in practice, even if we
use a specialized library. (Part of the reason is the frequent need to
reduce fractions $a/b$ by taking the $\gcd$ of $a$ and $b$.)
It is thus best to implement the algorithm -- in particular,
procedure \textsc{SumByLin} and its subroutines -- using only integer
arithmetic. Doing so also makes it easier to verify that the integers used
all fit in a certain range ($|n|<2^{127}$, say), and of course also
helps them fit in that range, in that we can simplify fractions before
we code: $(a/ b c)/(d/b f)$ (say) becomes $a f/ b d$, represented by the
pair of integers $(a f, b d)$.

{\em Square-roots and divisions.} On typical current 64-bit architectures,
a division takes as much time as several multiplications, and a square-root
takes roughly as much time as one or two divisions. (These are obviously
crude, general estimates.) Here, by ``taking a square-root'' of $x$
we mean computing the representable number closest to $\sqrt{x}$, or
the largest representable number no larger than $\sqrt{x}$, where
``representable'' means ``representable in extended precision'', that is,
as a number $2^e n$ with $|n|<2^{128}$ and $e \in [-(2^{14}-1),2^{14}-1] - 63$.

Incidentally, one should be extremely wary of using
hardware implementations of any floating-point operations other than
the four basic operations and the square-root; for instance, an implementation
of $\exp$
can give a result that is {\em not} the representable number closest to
$\exp(x)$ for given $x$. Fortunately, we do not need to use any floating-point
operations other than the square-root. The IEEE 754 standard requires that taking a square-root be implemented
correctly, that is, that the operation return the representable number closest
to $\sqrt{x}$, or the largest representable number $\leq \sqrt{x}$,
or the smallest such number $\geq \sqrt{x}$, depending on how we set the
rounding mode.

We actually need to compute $\lfloor \sqrt{n}\rfloor$ for $n$
a 128-bit integer. (We can assume that $n<2^{125}$, say.)
We do so by combining a single iteration of the procedure
in \cite{Zimmermann} (essentially Newton's method) with a
hardware implementation of a floating-point
extended-precision square-root in the sense we have just described.

It is of course in our interest to keep the number of divisions
(and square-roots) we perform
as low as possible; keeping the number of multiplications small is of course
also useful. Some easy modifications help: for instance, we can
conflate functions \textsc{Special1} and \textsc{Special0B} into a single
procedure; the value of $\gamma_1$ in the two functions differs by
exactly $m$.

{\em Parallelization.}
We parallelized the algorithm at two crucial places: one is function
\textsc{SArr} (Algorithm \ref{alg:sarr}), as we already discussed at the end of
\S \ref{sec:largenonfree}; the other one is function \textsc{DDSum}
(Algorithm \ref{alg:DDSum}), which involves a double loop. The task
inside the double loop (that is, \textsc{DoubleSum} or
\textsc{BruteDoubleSum}) is given to a processing element to compute on its
own. How exactly the double loop is traversed and parcelled out
is a matter that involves not just the usual trade-off between time and
space but also a possible trade-off between either and efficiency of
parallelization.

More specifically: it may be the case that the number of processing elements
is greater than the number of iterations of either loop
($\lceil (A'-A)/\Delta\rceil$ and $\lceil (B'-B)/\Delta\rceil$,
respectively), but smaller
than the number of iterations of the double loop. In that case,
parallelizing only the inside loop or the outside loop leads to an
under-utilization of processing elements. One alternative is
a na\"{i}ve parallelization of the double loop, with each processing
element recomputing the arrays $\mu$, $\mu'$ that it needs. That actually
turns out to be a workable solution: while recomputing arrays in this
way is wasteful, the overall time complexity does not change, and the total
space used is $O(\nu \Delta \log \log \max(A',B'))$, where $\nu$ is the
number of threads; this is slightly less space than $\nu$ instances
of \textsc{SumbyLin} use anyhow.

The alternative of computing and storing
the whole arrays $\mu$, $\mu'$ before entering the double loop
would allow us not to recompute them, but it would lead to
using (shared) memory on the order of $\max(A',B') \log \log \max(A',B')$,
which may be too large. Yet another alternative is to split the double
loop into squares of side about $\sqrt{\nu} \Delta$; then each array
segment $\mu$, $\mu'$ is recomputed only about $(A'-A)/(\sqrt{\nu} \Delta)$
or $(B'-B)/(\sqrt{\nu} \Delta)$ times, respectively, and we use
$O(\sqrt{\nu} \Delta)$ shared memory. Our implementation of this last
alternative, however, led to a significantly worse running time, at least
for $x=10^{19}$; in the end, we went with the ``workable solution'' above.
In the end, what is best may depend
on the parameter range and number of threads one is working with.

\section{Numerical results}

We computed $M(x)$ for $x=10^n$, $n\leq 23$, and
$x=2^n$, $n\leq 75$, beating the
records in \cite{Kuznetsov11} and \cite{Hurst18}. Our results
are the same as theirs, except that we obtain a sign
opposite to that in \cite[Table 1]{Kuznetsov11} for $x=10^{21}$; presumably
\cite{Kuznetsov11} contains a transcription mistake.

\begin{center}
  \begin{tabular}{l|l}
    $x$ & $M(x)$\\ 
    $10^{17}$ & $-21830254$\\ 
    $10^{18}$ & $-46758740$\\
    $10^{19}$ & $899990187$\\
    $10^{20}$ & $461113106$\\
    $10^{21}$ & $-3395895277$\\
    $10^{22}$ & $-2061910120$\\
    $10^{23}$ & $62467771689$
  \end{tabular}
  \begin{tabular}{l|l}
    $x$ & $M(x)$\\ 
    $2^{68}$ & $2092394726$\\
    $2^{69}$ & $-3748189801$\\
    $2^{70}$ & $9853266869$\\
    $2^{71}$ & $-12658250658$\\
    $2^{72}$ & $9558471405$\\
    $2^{73}$ & $-6524408924$\\
    $2^{74}$ & $-6336351930$\\
    $2^{75}$ & $-4000846218$
  \end{tabular}
\end{center}

Computing $M(x)$ for $x=10^{23}$ took about $18$ days and $14.6$ hours
on a 80-core machine
(Intel Xeon 6148, 2.40 GHz) shared
with other users. Computing $M(x)$ for $x=2^{75}=3.777\dotsc \cdot 10^{22}$
took about $9$ days and $16$ hours on the same machine.
As we shall see shortly, one parameter $c$ was more strictly
constrained for $x=10^{23}$, since
we needed to avoid overflow; 
we were able to optimize $c$ more freely for $2^{75}$.

\begin{figure}
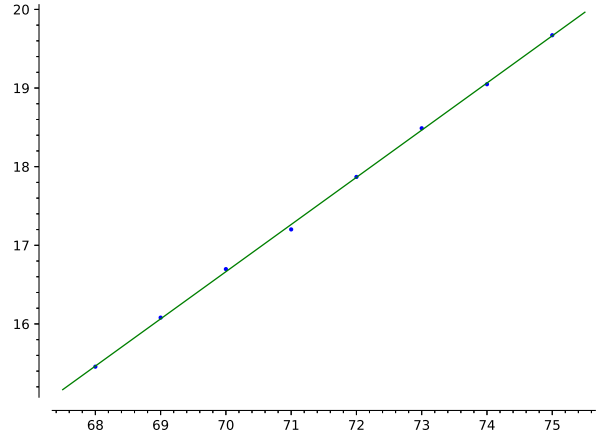
\label{fig:runtime}
  \begin{center}
    \sageplot[scale=.5]{plt}
  \end{center}
  \caption{Logarithm base $2$ of running time for input $x=2^n$}
\end{figure}

For a fixed choice of parameters, running time scaled approximately as $x^{3/5}$.
See Figure \ref{fig:runtime} for a plot\footnote{The first time we ran
  the program for $x=2^{75}$, we obtained a substantially higher running time,
  on the order of fourteen and a half days
  (as was reported on the first public draft of this paper).
  The time taken for $x=2^{71}$ was also higher on a first run, by about 20\%.
  We do not know the reason for this discrepancy, though demands by other users
  are probably the reason for $x=2^{71}$ and possibly also for $x=2^{75}$.
}
of the logarithm base $2$ of the running
time (in seconds) for $x= 2^n$, $n=68,69,\dotsc,75$
with $v= x^{2/5}/3$. We have drawn a line of slope $3/5$, with constant coefficient chosen by least squares to fit the points with
$68\leq n\leq 75$.

We also ran our code for $x=2^n$, $68\leq n\leq 75$, on a 128-core machine
based on two AMD EPYC 7702 (2GHz) processors. The results were of course the same as on the first computer, but running
time scaled more poorly,
particularly when passing from $2^{73}$ to $2^{74}$. (For whatever reason,
the program
gave up on $n=2^{75}$ on the second computer.)
The percentage of total time taken by the case of a large non-free variable was also much larger than on the first computer, and 
went up from $2^{73}$ to $2^{74}$. The reason for the difference in running times
in the two computers presumably lies in the differences between their respective
 memory architectures.
The dominance (in the second computer) of the case of a large non-free variable, whose usage of sieves is the most memory-intensive part of the program,
supports this diagnosis. It would then be advisable, for the sake
of reducing running times in practice, to improve on the memory usage of that
part of the program, either replacing \textsc{SegFactor} by the
improved sieve in \cite{Helfgott2020} -- sharply reducing
memory usage at the cost of increasing the asymptotic running time slightly,
as we have discussed -- or using a cache-efficient implementation of the
traditional segmented sieve as in \cite[Algorithm 1.2]{OSHP}. These two
strategies could be combined.


\begin{sagesilent}
  list = [(68,44959),(69,69284),(70,106258),(71,150765),(72,239472),(73,367780),
  (74,542183),(75,835136)]
  list2 = [(x[0],log(x[1])/log(2)) for x in list]
  bestt = gp(sum([(3/5)*x[0]-x[1] for x in list2])/len(list2))
  plt = points(list2)+plot(lambda x: (3/5)*x-bestt,67.5,75.5,color='green')
\end{sagesilent}

{\em Checking for overflow.}
Since our implementation uses 128-bit signed integers, it is crucial that
all integers used be of absolute value $<2^{127}$. What is critical here is
the quantity
\[\frac{\beta}{\delta} =
\frac{(\overline{x (m_\circ- (m-m_\circ))} /m_\circ^2 n_\circ  - r_0/q}{
- x/m_\circ n_\circ^2 - a/q} = 
\frac{(\overline{x (2 m_\circ - m)} q  - r_0 m_\circ^2 n_\circ) n_\circ }{(-x q - a m_\circ n_\circ^2) m_\circ}
\]
in \textsc{SumByLim}, where we write here $\overline{y}$ for the integer in
$\{0,1,\dotsc, m_\circ^2 n_{\circ}-1\}$ congruent to $y$
modulo $m_\circ^2 n_\circ$. The numerator could be as large as
$q m_\circ^2 n_\circ^2$
(The denominator is much smaller, since
$|-x/m_\circ n_\circ^2 - a/q|\leq 1/2 b q$.)
Since $q\leq 2 b$, $b\leq (A^4/6 x)^{1/3} \leq (v^4/6 x)^{1/3}$, $m_\circ,n_\circ
\leq v$ and $v = c x^{2/5} \frac{(\log \log x)^{3/5}}{(\log x)^{3/5}}$, we see that
\begin{equation}\label{eq:donrodrigo}
  q m_\circ^2 n_\circ^2\leq \frac{2 v^{16/3}}{(6 x)^{1/3}}
  = \frac{ 2 c^{16/3}}{6^{1/3}} \cdot x^{9/5}
\frac{(\log \log x)^{\frac{16}{5}}}{(\log x)^{\frac{16}{5}}}
.\end{equation}
\begin{sagesilent}
  cei=(lambda x,d : N((((RBF(10^d)*RBF(x)).upper()).ceil())/10^d))
  flo=(lambda x,d : N((((RBF(10^d)*RBF(x)).lower()).floor())/10^d))
  import re
  def pf(x):
    return sage.misc.latex.LatexExpr(re.sub(r'(\d)\.?0+($|\ )',r'\1\2',latex(RR(x))).replace(r'\times',r'\cdot'))
\end{sagesilent}
For $c=3/2$ 
and $x=2^{75} = 3.777\dotsc \cdot 10^{22}$,
\begin{sagesilent}
  f(c,x)= log((2*c^(16/3)/6^(1/3))*x^(9/5)*(log(log(x))/log(x))^(16/5))/log(2)
  lboo1 = flo(f(3/2,2^75),3)
  lboo2 = flo(f(9/8,10^23),3)
\end{sagesilent}
\[\log_2 \left(\frac{ 2 c^{16/3}}{6^{1/3}} x^{9/5}
\frac{(\log \log x)^{\frac{16}{5}}}{(\log x)^{\frac{16}{5}}}\right)
= \sage{pf(lboo1)}\dotsc < 127;\]
for $c = 9/8$ and $x=10^{23}$,
\[\log_2 \left(\frac{ 2 c^{16/3}}{6^{1/3}} x^{9/5}
\frac{(\log \log x)^{\frac{16}{5}}}{(\log x)^{\frac{16}{5}}}\right)
= \sage{pf(lboo2)}\dotsc < 127.\]
Thus, our implementation should give a correct result for $x = 10^{23}$,
for the choice $c=9/8$.
One can obviously go farther by using wider (or arbitrary-precision) integer
types.

There is another integer that might seem to be possibly larger, namely
the discriminant $\Delta = b^2-4 a c$ in the quadratic equations
solved in \textsc{QuadIneqZ}, which is called by functions
\textsc{Special1} and \textsc{Special0B}. However, that discriminant
is smaller than it looks at first.

The coefficient $\gamma_1$ in \textsc{Special0B} is
\[\begin{aligned}
&(-\lfloor R_0\rfloor q - r_0 + a_0 n_\circ) m =
(-\lfloor R_0\rfloor q - (\{R_0\} - \beta) q + a_0 n_\circ) m \\ &=
\left(- \left(\frac{x}{m_\circ n_\circ} - \frac{x}{m_\circ^2 n_\circ} (m-m_\circ)
\right)  q + \beta q + a_0 n_\circ\right) m\\
&=
\left(- \left(\frac{x}{m_\circ n_\circ} - \frac{x}{m_\circ^2 n_\circ} (m-m_\circ)
\right)  + \beta  +
\left(-\frac{x}{m_\circ n_\circ^2} - \delta\right) n_\circ\right) m q\\
&= \left(-\frac{2 x}{m_\circ n_\circ} +
\frac{x (m-m_\circ)}{m_\circ^2 n_\circ} +
O^*\left(\frac{1}{2 q}\right) + O^*\left(\frac{1}{2 b q}\right)
n_\circ\right) m q
.
\end{aligned}\]
Here the second term is negligible compared to the first one, and the
third term is negligible compared to the fourth one. We know that
\[\begin{aligned}
&\frac{x}{m_\circ n_\circ} m q \leq \frac{x}{m_\circ n_\circ} (m_\circ+a) \cdot 2 b
\leq \frac{2 b x}{n_\circ} + \frac{2 a b x}{m_\circ n_\circ}
\leq 2 x \sqrt[3]{\frac{A}{6 x}} + 2 x \sqrt[3]{\frac{A^2}{(6 x)^2}}
\\&\leq 2 x \sqrt[3]{\frac{v}{6 x}} + 2 x \sqrt[3]{\frac{v^2}{(6 x)^2}}
\leq 2 \sqrt[3]{\frac{c}{6}} \cdot x^{\frac{4}{5}}
\left(\frac{\log \log x}{\log x}\right)^{1/5} +
2 \left(\frac{c}{6}\right)^{\frac{2}{3}} x^{\frac{3}{5}}
\left(\frac{\log \log x}{\log x}\right)^{2/5}.
\end{aligned}\]
We also see that
\[\frac{n_\circ m}{2 b} \leq \frac{n_\circ m_\circ}{b}
\leq \sqrt[3]{6 x \cdot A^2} \leq \sqrt[3]{6 v^2 x}
\leq \sqrt[3]{6 c^2}\cdot x^{\frac{3}{5}}
\left(\frac{\log \log x}{\log x}\right)^{2/5}
.\]
The dominant term is thus $2 (c/6)^{1/3} x^{4/5} ((\log \log x)/\log x)^{1/5}$.
The coefficient $\gamma_1$ in \textsc{Special1} is equal to the
one we just considered, minus $m$, and thus has the same dominant term.

As for the term $- 4 a c$ (or $- 4 \gamma_0 \gamma_2$, so as not to conflict
with the other meanings of $a$ and $c$ here), it equals $4$ times
\[
a m x q = \frac{a}{q} m x q^2 =
\left(-\frac{x}{m_\circ n_\circ^2} - \delta\right) m x q^2
= - \frac{x^2 q^2 m}{m_\circ n_\circ^2} + O^*(m x).\]
Since
\[\frac{x^2 q^2}{n_\circ^2} \leq \frac{4 x^2 b^2}{B^2}
= 4 x^2 \sqrt[3]{A^2}{(6 x)^2}
\leq \frac{4}{6^{2/3}} x^{4/3} v^{2/3} \leq
\frac{4 c^{2/3}}{6^{2/3}} x^{8/5} \left(\frac{\log \log x}{\log x}\right)^{2/5}\]
and $m x\leq v x \leq c x^{7/5} (\log \log x)^{3/5}/(\log x)^{3/5}
$, we see that the main term here is at most
\[\frac{16 c^{2/3}}{6^{2/3}} x^{8/5}
\left(\frac{\log \log x}{\log x}\right)^{2/5}.\]

Since the two expressions we have just considered have opposite sign,
we conclude that 
the main term in the discriminant $\gamma_1^2 - 4 \gamma_0 \gamma_2$
is thus at most $(16 c^{2/3}/6^{2/3}) x^{8/5} (\log \log x)^{2/5}/(\log x)^{2/5}$,
that is, considerably
smaller than the term in \eqref{eq:donrodrigo}, at least for $x$
larger than a constant. For $c=3/2$ and $x=2^{75}$,
\begin{sagesilent}
  g(c,x)= log((16*c^(2/3)/6^(2/3))*x^(8/5)*(log(log(x))/log(x))^(2/5))/log(2)
  lbor1 = flo(g(3/2,2^75),3)
  lbor2 = flo(g(9/8,10^23),3)
\end{sagesilent}
\[\log_2 \frac{16 c^{2/3}}{6^{2/3}} x^{8/5}
\left(\frac{\log \log x}{\log x}\right)^{2/5}
= \sage{pf(lbor1)}\dotsc .\]
For $c=9/8$ and $x=10^{23}$,
\[\log_2 \frac{16 c^{2/3}}{6^{2/3}} x^{8/5}
\left(\frac{\log \log x}{\log x}\right)^{2/5}
= \sage{pf(lbor2)}\dotsc,\]
and thus we are out of danger of overflow for those parameters as well.

\appendix
\section{A sketch of an alternative algorithm}\label{subs:appalt}

As we mentioned in the introduction, we originally developed an algorithm
taking time $O(x^{3/5} (\log x)^{8/5})$ and
space $O(x^{3/10} \log x)$,
or, if the sieve in \cite{Helfgott2020} is used to factorize
integers in function \textsc{SArr} (Algorithm \ref{alg:sarr}),
time $O(x^{3/5} (\log x)^{8/5})$ and space
$O(x^{1/5} (\log x)^{1/5+5/3})$. The algorithm actually had an idea
in common with \cite{Helfgott2020}; as explained there, it
is an idea inspired by Vorono\"{i} and
Vinogradov's approach to the divisor problem.

Part of the improvement over that older algorithm resides in a better
(yet simple) procedure for computing sums of the form
$\sum_{d|n: d\leq a} \mu(d)$ (see Algorithm \ref{alg:factofun}); we analyzed it in \S \ref{sec:largenonfree}.
Other than that, the difference lies mainly in the computation
of the sum of $\mu(m) \mu(n) \lfloor x/m n\rfloor$ for
$(m,n)$ in a neighborhood $U = I_x\times I_y$
(see \S \ref{subs:handl} and
Algorithm \ref{alg:sumweight}).
Let us use the notation in
\S \ref{subs:handl}. In particular, write $I_x = [m_0-a,m_0+a)$,
  $I_y = [n_0-b,n_0+b)$. We have sums $S_0$, $S_1$, $S_2$, where
    $S_0$ is easy to compute and $S_2$ is the sum that we actually want
  to determine.

  In the version given in the current version of the paper, we
  compute the difference $S_1-S_0$ in time $O(a+b)$ and space $O(b \log b)$.
  Computing the difference $S_1-S_0$ in time $O((a+b) \log b)$ and space
  $O(b \log b)$
  (as we did in the previous version of the paper)
  is not actually hard; the main steps are: (i) sort
  the list of all pairs $(\{c_y (n-n_0)\},n)$ by their first element
  $\{c_y (n-n_0)\}$, (ii) use the sorted list to compute
  the sums $\sum_{n: \{c_y n\}\geq
    \{c_y n'\}} \mu(n)$ for different $n'$, and then (iii)
  search through the list
  as needed to determine the sum $\sum_{n:\{c_y n\}\geq \beta} \mu(n)$
  for any given value of $\beta$.

  The crux is how to compute $S_2-S_1$. In the current version, we analyze
  this difference with great care, after having determined
  the (at most) two arithmetic progressions in which the terms of
  $S_2-S_1$ that are non-zero must be contained. In the older version,
  we determined those arithmetic progressions in the same way as here
  (namely, by finding a Diophantine approximation $a/q$ to $c_y$). Within those progressions,
  however, we did not establish precisely what the non-zero terms were,
  but simply showed that they had to be contained in an interval $I\subset I_y$.
  We also showed that, for $q$ small, the interval $I$ had to be small as well,
  at least on average. (The number of elements
  of an arithmetic progression modulo $q$ within $I_y$ is $O(b/q)$, and so
  the case of $q$ large is not the main worry.)
  It is here that the argument in
  \cite[Ch.~III, exer.~3-6]{Vinogradov} came in handy:
  as we move from neighborhood to neighborhood, the quantity $c_y$
  keeps changing at a certain moderate speed, monotonically;
  thus, $c_y \mo \mathbb{Z}$ cannot spend too much time
  in major arcs on the circle $\mathbb{R}/\mathbb{Z}$.
  Only when
  $c_y \mo \mathbb{Z}$ lies in the major arcs can $q$ be small and the interval
  $I$ be large. Thus, just as claimed, the case of $q$ small and $I$ large occurs for
  few neighborhoods.

  We can thus simply determine $I$, and compute the terms that lie in the
  intersection of either of those two arithmetic progressions and
  their corresponding intervals $I$, and sum those terms. The time will be
  about $O(ab/q)$, unless $q$ is small, in which case one can do better,
  viz., $O(a |I|/q)$ or so. (Compare with the corresponding bound for the
  newer algorithm, namely, $O(a+b)$.) On average, we obtained savings of
  a factor of $O((\log b)/b)$, rather than $O(1/b)$, as we do now.

  Whether or not we use \cite{Helfgott2020} to factor integers $n\leq x/v$,
  we set $v = c x^{2/5}/(\log x)^{3/5}$, for $c$ a constant of our choice.
\pagebreak

\section{Pseudocode for algorithms}
In this section, we present the pseudocode for the algorithms referenced in this paper. To aid the reader, we begin with a diagram demonstrating the relationship between the algorithms.

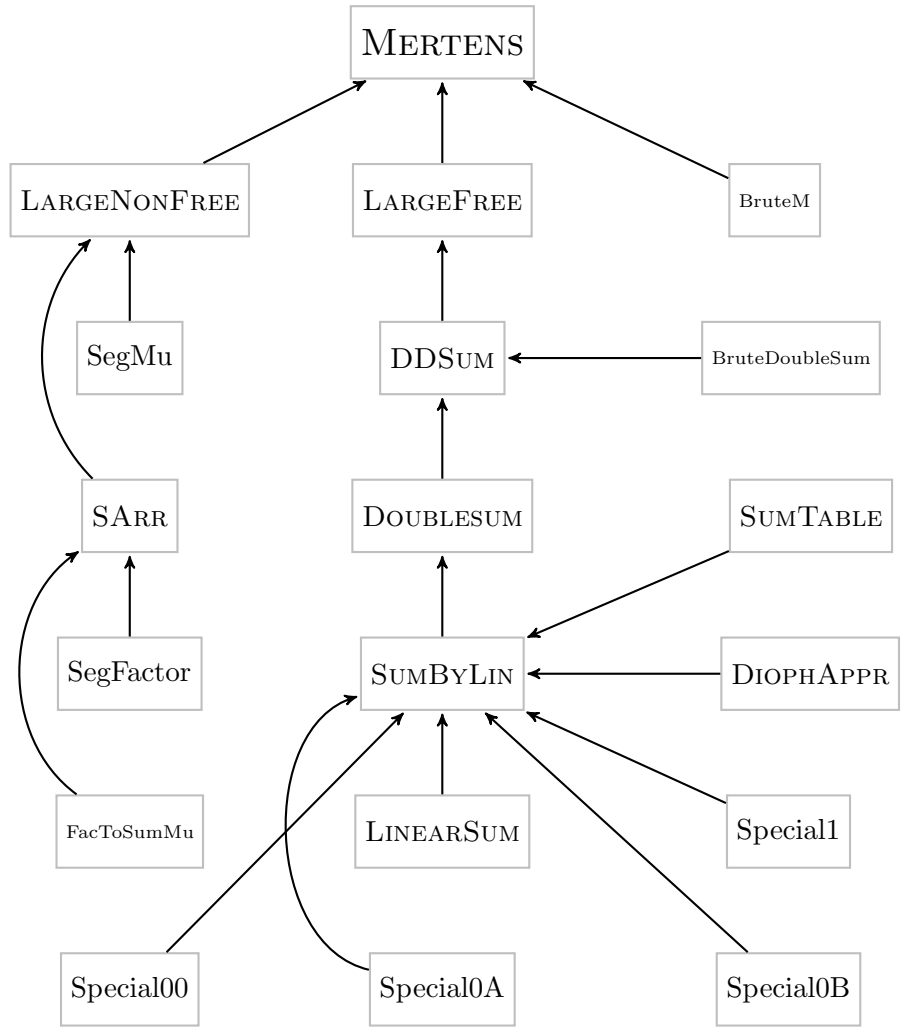
\begin{figure}[H]
\begin{tikzpicture}[auto,thick,node distance=2.5cm]

\tikzstyle{every state}=[rectangle,fill=none,draw=lightgray,text=black]

\begin{scope}[node distance=2.1cm and 2.6cm]
\node[state] (1)                    {\Large{\textsc{Mertens}}};
\node[state] (2) [below of=1]      {\textsc{LargeFree}};
\node[state] (4) [right=of 2]       {\tiny BruteM};
\node[state] (6) [below of=2]       {\textsc{DDSum}};
\node[state] (14) [left=of 6]    {SegMu};
\node[state] (3) [above of=14]        {\textsc{LargeNonFree}};
\node[state] (15) [below of=6]      {\textsc{Doublesum}};
\node[state] (8) [below of=15]       {\textsc{SumByLin}};
\node[state] (9) [right=of 8]       {\textsc{DiophAppr}};
\node[state] (7) [right=of 6]       {\tiny BruteDoubleSum};
\node[state] (10) [below of=8]      {\textsc{LinearSum}};
\node[state] (11) [right=of 15]     {\textsc{SumTable}};
\node[state] (12) [below of=14]     {\textsc{SArr}};
\node[state] (13) [below of=12]    {SegFactor};
\node[state] (5) [below of=13]       {\tiny FacToSumMu};
\node[state] (16) [right=of 10]     {Special1};
\node[state] (17) [below of=5]      {Special00};
\node[state] (18) [below of=10]     {Special0A};
\node[state] (19) [below of=16]     {Special0B};
\end{scope}

\begin{scope}[->,>=stealth',shorten >=1pt,]
  \path (2) edge  node {} (1);
  \path (3) edge  node {} (1);
  \path (4) edge  node {} (1);
  \path (5) edge [bend left=55] node {} (12);
  \path (6) edge node {} (2);
  \path (7) edge node {} (6);
  \path (15) edge node {} (6);
  \path (8) edge node {} (15);
  \path (9) edge node {} (8);
  \path (10) edge node {} (8);
  \path (11) edge node {} (8);
  \path (12) edge [bend left=45] node {} (3);
  \path (13) edge node {} (12);
  \path (14) edge node {} (3);
  \path (16) edge node {} (8);
  \path (17) edge node {} (8);
  \path (18) edge [bend left=75] node {} (8);
  \path (19) edge node {} (8);
  \end{scope}

\end{tikzpicture}
\caption{Dependency diagram}
\end{figure}

\begin{algorithm}
  \caption{Main algorithm: compute $M(x) = \sum_{n\leq x} \mu(n)$}
  \label{alg:mertens}
  \begin{algorithmic}[1]
    \Function{Mertens}{$x$}
    \Ensure{$\sum_{n\leq x} \mu(n)$}
    \State{$c\gets 3/2$} \Comment{hand-tuned value, change at will}
    \State{$u= \sqrt{x}$, $v\gets c x^{2/5} (\log \log x)^{3/5}/(\log x)^{3/5}$}
    \State{$M \gets 2\cdot \textsc{BruteM}(u)$}
    \State{$M\gets M -
        \textsc{LargeNonFree}(x,v,u) -
        \textsc{LargeFree}(x,v)$}
    \State{\Return{$M$}}
    \EndFunction
    \vskip 3pt
    \noindent {\bf Time:} $O\left(x^{\frac{3}{5}} (\log x)^{3/5}
    (\log \log x)^{2/5}\right)$.
    \vskip 1pt
    \noindent      {\bf Space:}
    $O\left(x^{\frac{3}{10}} (\log x)^{\frac{13}{10}} (\log \log x)^{-\frac{3}{10}}\right)$.
      \end{algorithmic}
\end{algorithm}

\begin{algorithm}
  \caption{Compute $M(x) = \sum_{n\leq x} \mu(n)$ by brute force}
  \label{alg:brutemert}
  \begin{algorithmic}[1]
    \Function{BruteM}{$x$}
    \Ensure{$\sum_{n\leq x} \mu(n)$}
    \State{$M\gets 0$, $\Delta\gets \lfloor \sqrt{x}\rfloor$}
    \For{$0\leq j< \lceil x/\Delta\rceil$}
    \State{$n_0\gets j \Delta + 1$}
    \State{$\mu\gets \textsc{SegMu}(n_0,\Delta)$}
    \For{$n_0\leq n\leq \min(n_0+\Delta-1,x)$}
    \State{$M\gets M + \mu_{n-n_0}$}
    \EndFor
    \EndFor
    \State{\Return{$M$}}
    \EndFunction
    \vskip 3pt
    \noindent {\bf Time:} $O(x \log \log x)$.\;
              {\bf Space:} $O(\sqrt{x} \log x)$.
              \vskip 1pt \noindent
      \end{algorithmic}
\end{algorithm}

\begin{algorithm}
  \caption{The case of a large non-free variable}\label{alg:largvar}
  \begin{algorithmic}[1]
    \Function{LargeNonFree}{$x$,$v$,$u$}
    \Ensure{
      $\sum_{n\leq x} \sum_{m_1 m_2 n_1 = n: m_1, m_2\leq u, \max(
   m_1,m_2)>v} \mu(m_1) \mu(m_2)$}
    \State{$n_0\gets \lfloor u\rfloor + 1$,
      $r_0\gets \lfloor x/(\lfloor u\rfloor + 1)\rfloor + 1$}
    \State{$\Delta\gets \lceil \sqrt{\max(u,x/v)}\rceil$,
          $\mathbf{S}\gets \textsc{SArr}(x,r_0,\Delta,1)$}
    \State{$\Sigma\gets 0$, $\sigma\gets 0$}
    \For{$n = \lfloor u\rfloor, \lfloor u\rfloor-1,\dotsc,\lfloor v\rfloor + 1$}
    \If{$n< n_0$}
    \State{$n_0\gets \max(n_0-(\Delta+1),1)$, $\mu\gets \textsc{SegMu}(n_0,\Delta)$}
    \EndIf
    \State{$\sigma\gets \sigma+\mu_{n-n_0} \lfloor x/n^2\rfloor$}
    \While{$x/n>r_0+\Delta$}
    \State{$r_0\gets r_0+\Delta+1$, $\mathbf{S}\gets \textsc{SArr}(x,r_0,\Delta,\mathbf{S}_{\Delta})$}
    \EndWhile
    \State{$\Sigma\gets \Sigma + 2 \mu_{n-n_0} \cdot \left( - \sigma + \mathbf{S}_{\left\lfloor \frac{x}{n}\right\rfloor - r_0}\right) + \mu_{n-n_0}^2 \left\lfloor
      x/n^2\right\rfloor$}
    \EndFor
    \State{\Return{$\Sigma$}}
    \EndFunction
    \vskip 3pt
    \noindent {\bf Time:} $O\left(\left(\frac{x}{v} + u\right) \log \log x
    \right)$ 
    \vskip 1pt
    \noindent          {\bf Space:} $O\left(\sqrt{\max(x/v,u)}\cdot \log x\right)$.
  \end{algorithmic}
\end{algorithm}

\begin{algorithm}
  \caption{Compute the main sum needed for \textsc{LargeNonFree}}\label{alg:sarr}
  \begin{algorithmic}[1]
    \Function{SArr}{$x$,$r_0$,$\Delta$,$S_0$}
    \Ensure{for $0\leq j\leq \Delta$,
      $\mathbf{S}_j = \sum_{r\leq r_0+j} \sum_{b|r: b\leq \frac{x}{r}} \mu(b)$.}
    \Require{$\mathbf{S}_0= \sum_{r<r_0} \sum_{b|r:b\leq \frac{x}{r}} \mu(b)$}
    \State{$F\gets \textsc{SegFactor}(r_0,\Delta)$, $S\gets S_0$}
    \For{$r=r_0,r_0+1,\dotsc,r_0+\Delta$}
    \State{$S\gets S + \textsc{FacToSumMu}(\mathbf{F}_{r-r_0},x/r)$,
      $\mathbf{S}_{r-r_0}\gets S$}
    \EndFor
    \State{\Return{$\mathbf{S}$}}
        \EndFunction
    \vskip 3pt
    \noindent {\bf Time:} $O\left((\sqrt{r_0} +\Delta) \log \log x
    \right)$ 
    \vskip 1pt
    \noindent          {\bf Space:} $O\left((\sqrt{r_0}+\Delta)\log x\right)$.
         \end{algorithmic}
\end{algorithm}

\begin{algorithm}
  \caption{The case of a large free variable}\label{alg:smallnonfree}
  \begin{algorithmic}[1]
    \Function{LargeFree}{$x$,$v$}
        \Ensure{
      $\sum_{n\leq x} \sum_{m_1 m_2 n_1 = n:\; m_1, m_2\leq v} \mu(m_1) \mu(m_2)$}
        \State{$S\gets 0$, $A'\gets \lfloor v\rfloor+1$, $C\gets 10$,
        $D\gets 8$}\Comment{$C$ and $D$ are hand-tuned}
    \While{$A'\geq \max(2 (6 C^3 x)^{1/4},\lceil \sqrt{v}\rceil,2 D)$}
    \State{$B'\gets A'$, $A\gets A' - 2 \lfloor A'/2 D\rfloor$}
    \While{$B'\geq \max(2 (6 C^3 x/A)^{1/3},\lceil \sqrt{v}\rceil,2 D)$}
    \State{$B\gets B' - 2 \lfloor B'/2 D\rfloor$}
    \State{$a\gets \sqrt[3]{\frac{A^4}{6 x}}$,
      $b\gets \sqrt[3]{\frac{A B^3}{6 x}}$,
      $\Delta \gets
\lceil \sqrt{v}/\max(2 a, 2 b)\rceil\cdot \max(2 a, 2 b)$}
    \State{$S\gets S + \textsc{DDSum}(A,A',
      B,B',x,\Delta,1,a,b)\cdot \begin{cases} 1 &\text{if $A=B$,}\\
        2 &\text{if $A>B$.}\end{cases}$}
    \State{$B'\gets B$}
    \EndWhile
    \State{$S\gets S+2\cdot \textsc{DDSum}(A,A',
      1,B',x,\lceil \sqrt{v}\rceil,0,0,0)$}
   \State{$A'\gets A$}
    \EndWhile
    \State{$S\gets S + \textsc{DDSum}(A,A',
      1,B',x,\lceil \sqrt{v}\rceil,0,0,0)$}
    \State{\Return{$S$}}
    \EndFunction
       \vskip 3pt
       \noindent {\bf Time:} $O\left(v^{2/3} x^{1/3} \log x + v^{3/2} \log \log x\right)$
             \vskip 1pt
             \noindent {\bf Space:} $O\left(\sqrt{v} \log \log x
             + (v^4/x)^{1/3} \log x\right)$
    \end{algorithmic}
\end{algorithm}

\begin{algorithm}\label{alg:DDSum}
  \caption{split
    $\sum_{(m,n) \in [A,A')\times [B,B')}
        \mu(m) \mu(n) \left\lfloor \frac{x}{m n}\right\rfloor
        $ into smaller sums}
  \label{alg:ddsum}
  \begin{algorithmic}[1]
    \Function{DDSum}{$A$,$A'$,$B$,$B'$,$x$,$\Delta$,$\gamma$,$a$,$b$}
    \Ensure{$\sum_{(m,n) \in [A,A')\times [B,B')} \mu(m) \mu(n)
           \left\lfloor \frac{x}{m n}\right\rfloor$}
    \Require{$A,B\geq 1$, $2|\Delta$, $A'\equiv A \mo 2$, $B'\equiv B \mo 2$}
    \State{$S\gets 0$}
    \For{$m_0 \in \lbrack A,A') \cap (A + \Delta \mathbb{Z})$}
    \State{$m_1\gets \min(m_0+\Delta,A')$,
      $\mu\gets \textsc{SegMu}(m_0,\Delta)$}
    \For{$n_0 \in \lbrack B,B') \cap (B + \Delta \mathbb{Z})$}
    \State{$n_1\gets \min(n_0+\Delta,B')$,
      $\mu'\gets \textsc{SegMu}(n_0,\Delta)$}
    \If{$\gamma=1$}
    \State{$S \gets S + \textsc{DoubleSum}(m_0,m_1,n_0,n_1,a,b,\mu,\mu',x)$}
    \Else
    \State{$F(m,n):=\lfloor x/m n\rfloor$, $f(m) := \mu_{m-m_0}$, $g(n) := \mu_{n-n_0}'$}
    \State{$S \gets S + \textsc{BruteDoubleSum}(m_0,m_1,
        n_0,n_1,\mu,\mu',F)$}
    \EndIf
    \EndFor
    \EndFor
    \State{\Return{$S$}}
    \EndFunction
    \vskip 3pt
    \noindent {\bf Time:}
    $O\left(\left\lceil \frac{A'-A}{\Delta} \right\rceil
        \left\lceil \frac{B'-B}{\Delta} \right\rceil \Delta \log \log \Delta\right)$,
        assuming $\Delta\gg \sqrt{\max(A',B')}$, plus time taken by
        \textsc{DoubleSum} or \textsc{BruteDoubleSum}.
    \vskip 1pt \noindent
           {\bf Space:} $O(\Delta \log \log \max(A',B'))$,
           mainly from \textsc{SegMu}
      \end{algorithmic}
\end{algorithm}

\begin{algorithm}
  \caption{
    $\sum_{(m,n) \in [m_0,m_1)\times [n_0,n_1)}
                f(m) g(n) F(m,n)$ by brute force}
  \label{alg:brutesum}
  \begin{algorithmic}[1]
    \Function{BruteDoubleSum}{$m_0$,$m_1$,$n_0$,$n_1$,$f$,$g$,$x$}
    \Ensure{$\sum_{(m,n) \in [m_0,m_1)\times [n_0,n_1)} f(m) g(n)
      F(m,n)$}
     \State{$S\gets 0$}
    \For{$m_0\leq m< m_1$}
    \For{$n_0\leq n< n_1$}
    \State{$S \gets S + f(m) g(n) F(m,n)$}
    \EndFor
    \EndFor
    \State{\Return{$S$}}
    \EndFunction
    \vskip 3pt
    \noindent {\bf Time:} $O((m_1-m_0) (n_1-n_0) + 1)$.
    \vskip 1pt \noindent
              {\bf Space:} that of the inputs, plus $O(1)$.
      \end{algorithmic}
\end{algorithm}

\begin{algorithm}
  \caption{compute
    $\sum_{(m,n) \in [m_0,m_1)\times [n_0,n_1)}
    f_{m-m_0} g_{n-n_0} \left\lfloor \frac{x}{m n}\right\rfloor$
    }\label{alg:doublesum}
 \begin{algorithmic}[1]
    \Function{Doublesum}{$m_0$,$m_1$,$n_0$,$n_1$,$a$,$b$,$f$,$g$,$x$}
    \Ensure{
$\sum_{(m,n) \in [m_0,m_1)\times [n_0,n_1)}
    f_{m-m_0} g_{n-n_0} \left\lfloor \frac{x}{m n}\right\rfloor$}
  \Require{$m_0,n_0\geq 1$, $m_1\leq 2 m_0$, $n_1\leq 2 n_0$,
    $2|m_1-m_0$, $2|n_1-n_0$,
          and all conditions for
          \textsc{SumByLin}}
        \State{$S\gets 0$}
        \For{$0\leq j< \lceil (m_1-m_0)/2a\rceil$}
        \State{$m_-\gets m_0+j\cdot 2 a$,
          $m_+\gets \min(m_0 + (j+1)\cdot 2 a, m_1)$}
        \State{$m_\circ\gets (m_-+m_+)/2$,
          $m_\Delta\gets (m_+-m_-)/2$}\Comment{midpoint, width}
        \For{$0\leq k< \lceil (n_1-n_0)/2b\rceil$}
        \State{$n_-\gets n_0+k\cdot 2 b$,
          $n_+\gets \min(n_0+(k+1)\cdot 2 b, n_1)$}
        \State{$n_\circ\gets (n_-+n_+)/2$,
          $n_\Delta\gets (n_+-n_-)/2$}\Comment{midpoint, width}
        \State{$f(m):=f_{m+m_\circ-m_0}$, $g(n):=g_{n+n_\circ-n_0}$}
        \State{$S\gets S+\textsc{SumByLin}(f,g,x,m_\circ,n_\circ,a,b)$}
        \EndFor
        \EndFor
        \State{\Return{$S$}}
        \EndFunction
       \vskip 3pt
       \noindent {\bf Time:} $O\left(\frac{A B}{\min(a,b)}\right)$ 
      \vskip 1pt
      \noindent {\bf Space:} that of the inputs, plus $O(b\log b)$
    \end{algorithmic}
\end{algorithm}

\begin{algorithm}
  \caption{Finding a Diophantine approximation via continued fractions}\label{alg:diophfrac}
  \begin{algorithmic}[1]
    \Function{DiophAppr}{$\alpha$,$Q$}
    \Ensure{$(a,a^{-1},q,s)$ s.t.
      $\left|\alpha-\frac{a}{q}\right|\leq \frac{1}{q Q}$, $(a,q)=1$, $q\leq Q$, $a a^{-1} \equiv 1 \mo q$ and $s=\sgn(\alpha-a/q)$}
    \State{$b\gets \lfloor \alpha\rfloor$, $p\gets b$, $q\gets 1$,
    $p_-\gets 1$, $q_-\gets 0$, $s\gets 1$}
    \While{$q\leq Q$}
    \If{$\alpha=b$}
    \Return{$(p,- s q_-, q,0)$}
    \EndIf
    \State{$\alpha \gets 1/(\alpha-b)$}
    \State{$b\gets \lfloor \alpha\rfloor$,
      $(p_+,q_+) \gets b\cdot (p,q) + (p_-,q_-)$}
    \State{$(p_-,q_-)\gets (p,q)$, $(p,q)\gets (p_+,q_+)$, $s\gets -s$}
    \EndWhile
    \State{\Return{$(p_-, s q,q_-,-s)$}}
    \EndFunction
    \vskip 3pt
\noindent {\bf Time:} $O(\log \max(Q,\den(\alpha))$.\; {\bf Space:} $O(1)$.
      \end{algorithmic}
\end{algorithm}

\begin{algorithm}
  \caption{Preparing tables of partial sums by congruence class}\label{alg:sumtable}
  \begin{algorithmic}[1]
    \Function{SumTable}{$f$,$b$,$a_0$,$q$}
    \Ensure{$(F,\rho,\sigma)$ where $F_{n_0} =
     \sum_{-b\leq n\leq n_0: n\equiv n_0 \mo q} f(n)$ for
     $-b\leq n_0<b$}
    \Ensure{$\rho_r = \sum_{-b\leq n< b: a_0 n \equiv r \mo q} f(n)$ and
      $\sigma_r = \sum_{j=q-r+1}^{q-1} \rho_j$.}
    \Require{$q\leq 2 b$}    
    \For{$n\in [-b,-b+q)$}
      \State{$F_n\gets f(n)$}
      \EndFor
      \For{$n\in [-b+q,b)$}
        \State{$F_n\gets F_{n-q}+f(n)$}
        \EndFor
        \State{$r\gets \textsc{Mod}(a_0 (b-q),q)$}
        \For{$n\in \{b-q,\dotsc,b-1\}$}
        \State{$\rho_r \gets F_n$}
        \State{$r\gets \textsc{Mod}(r + a_0,q)$}
        \EndFor
        \State{$\sigma_0\gets 0$, $\sigma_1\gets 0$}
        \For{$r\in \{1,2,\dotsc,q-1\}$}
        \State{$\sigma_{r+1} \gets \sigma_r + \rho_{q-r}$}
        \EndFor
        \State{\Return{$(F,\rho,\sigma)$}}
        \EndFunction
  
\vskip 3pt
\noindent {\bf Time:} $O(b)$.\; {\bf Space:} $O(b\log b)$.

\vskip 5pt
\Function{RaySum}{$f$,$q$,$b$,$\delta$}
\State{$S\gets 0$}
\If{$\delta<0$}
\For{$n\in \left\{q, 2 q,\dotsc,\left\lfloor(b-1)/q\right\rfloor q\right\}$}
\State{$S\gets S + f[n]$}
\EndFor
\EndIf
\If{$\delta>0$}
\For{$n\in \{q, 2 q,\dotsc,\left\lfloor b/q\right\rfloor q\}$}
\State{$S\gets S + f[-n]$}
\EndFor
\EndIf
\State{\Return{$S$}}
\EndFunction

\vskip 3pt
      \noindent {\bf Time:} $O(n/q)$\;\;\;{\bf Space:} $O(1)$
      \vskip 5pt
      \Function{Mod}{$a$,$q$}
      \EndFunction
      
      Returns the integer $0\leq r<q$ such that $r\equiv a \mo q$.
           \vskip 2pt
      \noindent {\bf Time and space:} $O(1)$.
      \vskip 5pt

            \vskip 5pt
            \Function{Sgn}{$\delta$}
            \If{$\delta<0$}
            \State{\Return{$-1$}}
            \ElsIf{$\delta>0$}
            \State{\Return{$1$}}
            \Else
            \State{\Return{$0$}}
            \EndIf
      \EndFunction
      
      Returns the integer $0\leq r<q$ such that $r\equiv a \mo q$.
           \vskip 2pt
      \noindent {\bf Time and space:} $O(1)$.
      \vskip 5pt

\end{algorithmic}
\end{algorithm}

\begin{algorithm}
  \caption{Summing with a weight $x/m n$ using a linear approximation}\label{alg:sumweight}
  \begin{algorithmic}[1]

    \Function{SumByLin}{$f$,$g$,$x$,$m_\circ$,$n_\circ$,$a$,$b$}
    \Ensure{
      $\sum_{(m,n)\in U} f(m) g(n) \left\lfloor \frac{x}{(m+m_\circ)
        (n+n_\circ)}\right\rfloor$
      for $U = \lbrack -a,a)\times \lbrack -b,b)$,
      $a,b\in \mathbb{Z}^+$
    }
    \Require{the difference between $\frac{x}{(m+m_\circ) (n+n_\circ)}$ and its linear approximation
      around $(0,0)$ has absolute value $\leq 1/2 b$ on $U$}
    \State{$\alpha_0 \gets \frac{x}{m_\circ n_\circ}$,
      $\alpha_1 = - \frac{x}{m_\circ^2 n_\circ}$,
      $\alpha_2 = - \frac{x}{m_\circ n_\circ^2}$}
    \State{$S\gets \textsc{LinearSum}(f,g,a,b,\alpha_0,\alpha_1,\alpha_2)$}    
    \State{$(a_0,\overline{a_0},q, s)\gets \textsc{DiophAppr}(\alpha_2,2 b)$, $\delta \gets \alpha_2 - a_0/q$, $\delta' \gets \textsc{Sgn}(\delta)$}
    \State{$Z\gets \textsc{RaySum}(g,q,b,s_\delta)$}
    \State{$(G,\rho,\sigma)\gets \textsc{SumTable}(g,b,a_0,q)$}
    \For{$m\in \lbrack -a,a)$ such that $f(m)\ne 0$}
    \State{$R_0\gets \alpha_0 + \alpha_1 m$,
      $r_0\gets \lfloor \{R_0\} q +1/2\rfloor$,
     $m'\gets m_\circ+m$}
    \State{$\beta\gets \{R_0\} - r_0/q$, $\beta'\gets \textsc{Sgn}(\beta)$}
    \If{$\delta\ne 0$}
    \State{$Q\gets \beta/\delta$}\Comment{the value of $Q$ for $\delta=0$
      is arbitrary}
    \EndIf
    \State{$T \gets \sigma_{r_0} + \textsc{Special0A}(G,q,a_0,\overline{a_0},r_0,b,Q,\beta',\delta')$}
    \If{$q>1$}
    \State{$T\gets T + \textsc{Special1}(G,x,q,a_0,\overline{a_0},R_0,r_0,n_\circ,m',b)$}
    \State{$T\gets T + \textsc{Special0B}(G,x,q,a_0,\overline{a_0},R_0,r_0,n_\circ,m',b,Q,\beta',\delta')$}
    \Else
    \State{$T\gets T + \textsc{Special00}(G,x,q,a_0,\overline{a_0},R_0,r_0,n_\circ,m',b,Q,\delta')$}
    \EndIf
    \If{$0<r_0<q$}
    \State{$T \gets T +  Z$}
    \EndIf
    \State{$S \gets S + f(m)\cdot T$}
     \EndFor
    \State{\Return{$S$}}
\EndFunction
       \vskip 3pt
       \noindent {\bf Time:} $O(a+b)$
      \vskip 1pt
      \noindent {\bf Space:} $O(b \log b)$, mainly from \textsc{SumTable}
  \end{algorithmic}
\end{algorithm}

\begin{algorithm}
  \caption{Table lookup}\label{alg:tablook}
  \begin{algorithmic}[1]
    \Function{SumInter}{$G$,$r$,$I$,$b$,$q$}
    \Require{$I=[I_0,I_1]$, where $I_0,I_1\in \mathbb{Z}$, $I_0\leq I_1$,
      or $I=\emptyset$}
    \If{$I\ne \emptyset$}
    \State{\Return{$0$}}
    \EndIf
    \State{$r_0\gets \textsc{FlCong}(I_0-1,r,q)$, $r_1\gets \textsc{FlCong}(\min(I_1,b-1),r,q)$}
    \If{$(r_0>r_1) \vee (r_1<-b)$}
    \State{\Return{$0$}}
    \EndIf
    \If{$r_0\geq -b$}
    \State{\Return{$G_{r_1} - G_{r_0}$}}
    \Else
    \State{\Return{$G_{r_1}$}}
    \EndIf
\EndFunction
    \vskip 3pt
    \noindent {\bf Time and space:} $O(1)$.
    \vskip 5pt
        \Function{FlCong}{$n$,$a$,$q$}
    \Ensure{Returns largest integer $\leq n$ congruent to $a \mo q$}
        \State{\Return{$n-\textsc{Mod}(n-a,q)$}}
    \EndFunction
        \vskip 3pt
    \noindent {\bf Time and space:} $O(1)$.
  \end{algorithmic}
\end{algorithm}

\begin{algorithm}
  \caption{$L_2-L_1$ for special moduli: quadratic equations}\label{alg:specquad}
  \begin{algorithmic}[1]
\Function{Special1}{$G$,$x$,$q$,$a$,$\overline{a}$,$R_0$,$r_0$,$n_\circ$,$m$,$b$}
\State{$\gamma_1 = (-\lfloor R_0\rfloor q - (r_0+1) + a n_\circ) m$}
\State{$r\gets (-1-r_0) \overline{a}$}
\State{$I\gets \textsc{QuadIneqZ}(-a m,\gamma_1,x q)-n_0$}
\State{\Return{$\textsc{SumInter}(G,r,(-\infty,\infty),b,q)-
    \textsc{SumInter}(G,r,I,b,q)$}}
\EndFunction
\vskip 5pt
\Function{Special0b}{$G$,$x$,$q$,$a$,$\overline{a}$,$R_0$,$r_0$,$n_\circ$,$m$,$b$,$Q$,$s_\beta$,$s_\delta$}
\State{$\gamma_1 = (-\lfloor R_0\rfloor q - r_0 + a n_\circ) m$}
\State{$I\gets \textsc{QuadIneqZ}(-a m,\gamma_1,x q)-n_\circ$}
\If{$s_\delta>0$}
\State{$J \gets (-\infty,-\lfloor Q \rfloor - 1]$}
\ElsIf{$s_\delta<0$}
\State{$J \gets [-\lceil Q\rceil+1,\infty)$}
  \ElsIf{$s_\beta\geq 0$}
  \State{$J\gets \emptyset$}
  \Else
  \State{$J \gets (-\infty,\infty)$}
    \EndIf
\State{\Return{$\textsc{SumInter}(G,- r_0 \overline{a},J,b,q) - 
    \textsc{SumInter}(G,- r_0 \overline{a},I\cap J,b,q)$}}
\EndFunction
    \vskip 3pt
    \noindent {\bf Time and space:} $O(1)$.
  \end{algorithmic}
\end{algorithm}

\begin{algorithm}
  \caption{$L_2-L_1$: the case $q=1$}\label{alg:caseq1}
  \begin{algorithmic}[1]
\Function{Special00}{$G$,$x$,$q$,$a$,$\overline{a}$,$R_0$,$r_0$,$n_\circ$,$m$,$b$,$Q$,$s_\delta$}
    \If{$s_\delta>0$}
\State{$J \gets (-\infty,-\lfloor Q \rfloor - 1]$}
\ElsIf{$s_\delta<0$}
\State{$J \gets [-\lceil Q\rceil+1,\infty)$}
    \Else
    \State{$J\gets \emptyset$}
    \EndIf
    \For{$j=0,1$}
\If{$a\ne 0$}
    \State{$\gamma_1 = (-\lfloor R_0\rfloor - (r_0+j) + a n_\circ) m$}
    \State{$I_j\gets \textsc{QuadIneqZ}(-a m,\gamma_1,x)-n_\circ$}
    \Else
    \State{$I_j\gets (-\infty,\lfloor (x/m)/(\lfloor R_0\rfloor + r_0+j)\rfloor-n_\circ]$}
    \EndIf
    \EndFor
    \State{$S\gets
      \textsc{SumInter}(G,0,I_0\cap J,b,q)$}
    \State{$S\gets 
      S+\textsc{SumInter}(G,0,I_1\cap (\mathbb{R}\setminus J),b,q)$}
    \State{\Return{
$\textsc{SumInter}(G,0,(-\infty,\infty),b,q)-
        S$}}
    \EndFunction
    \vskip 3pt
    \noindent {\bf Time and space:} $O(1)$.
  \end{algorithmic}
\end{algorithm}
    
\begin{algorithm}
  \caption{$L_1-L_0$: casework for
  $a_0 (n-n_0) + r_0 \equiv 0 \mo q$}\label{alg:casewo}
\begin{algorithmic}[1]
  \Function{Special0a}{$G$,$q$,$a$,$\overline{a}$,$r_0$,$b$,$Q$,
    $s_\beta$,$s_\delta$}
\If{$0<r_0<q$}
\If{$s_\delta\ne 0$}
\If{$s_\delta> 0$}
\State{$I\gets [-\lfloor Q \rfloor,\infty)$}
\Else
\State{$I\gets (-\infty,-\lceil Q\rceil]$}
\EndIf
\ElsIf{$s_\beta\geq 0$}
\State{$I\gets (-\infty,\infty)$}
\Else
\State{$I\gets \emptyset$}
\EndIf
\Else
\If{$s_\delta=0 \vee s_\beta=0$}
\State{$I\gets \emptyset$}
\ElsIf{$s_\beta<0$}
\If{$s_\delta<0$}
\State{$S\gets \textsc{SumInter}(G,- r_0 \overline{a},
  (-\infty,-\lceil Q\rceil],b,q)$}
\State{\Return{$S+\textsc{SumInter}(G,- r_0 \overline{a},
           (0,\infty),b,q)$}}
\Else
\State{$S\gets \textsc{SumInter}(G,- r_0 \overline{a},(-\infty,0),b,q)$}
\State{\Return{$S+\textsc{SumInter}(G,- r_0 \overline{a},
           [-\lfloor Q \rfloor,\infty),b,q)$}}
\EndIf
\Else
\If{$s_\delta>0$}
\State{$I\gets [- \lfloor Q\rfloor,0)$}
  \Else
  \State{$I\gets (0,- \lceil Q \rceil]$}
\EndIf
\EndIf
\EndIf

\State{\Return{\textsc{SumInter}($G$,$- r_0 \overline{a}$,$I$,$b$,$q$)}}
\EndFunction

    \vskip 3pt
    \noindent {\bf Time and space:} $O(1)$.
  \end{algorithmic}
\end{algorithm}

\begin{algorithm}
  \caption{Summing with floors of linear expressions as weights}\label{alg:linearsum}
  \begin{algorithmic}[1]
    \Function{LinearSum}{$f$,$g$,$a$,$b$,$\alpha_0$,$\alpha_1$,$\alpha_2$}
\Ensure{$\sum_{(m,n)\in U} f(m) g(n) 
  (\lfloor \alpha_0 + \alpha_1 m\rfloor
  + \lfloor\alpha_2 n\rfloor)$ for $U = \lbrack -a,a)\times \lbrack -b,b)$}
\State{$S_1\gets 0$, $S_{1,0}\gets 0$, $S_2\gets 0$, $S_{2,0}\gets 0$}
\For{$m\in [-a,a)\cap \mathbb{Z}$}
\State{$S_1\gets S_1+ f[m]\cdot \lfloor \alpha_0 + \alpha_1 m\rfloor$,
$S_{1,0}\gets S_{1,0} + f[m]$}
\EndFor
\For{$n\in [-b,b)\cap \mathbb{Z}$}
\State{$S_2\gets S_2 + g[n]\cdot \lfloor\alpha_2 n\rfloor$, $S_{2,0}\gets S_{2,0} + g[m]$}
\EndFor
\State{\Return{$S_1\cdot S_{2,0} + S_{1,0}\cdot S_2$}}
\EndFunction
       \vskip 3pt
      \noindent {\bf Time:} $O(\max(a+1,b+1))$.
      \vskip 1pt
      \noindent {\bf Space:} $O(1)$ 
      \vskip 5pt
      
\end{algorithmic}
\end{algorithm}


\begin{algorithm}
  \caption{A little Babylonian routine}\label{alg:quadrineq}
    \begin{algorithmic}[1]
    \Function{QuadIneqZ}{$a$,$b$,$c$}
    \Ensure{Returns an interval $I$ such that}
    \Ensure{
      $I\cap \mathbb{Z} = \{x\in \mathbb{Z}: a x^2 + b x + c\geq 0\}$,
      if $a<0$,}
    \Ensure{
      $I\cap \mathbb{Z} = \{x\in \mathbb{Z}: a x^2 + b x + c< 0\}$,
      if $a>0$.}    
    \Require{$a,b,c\in \mathbb{Z}$, $a\ne 0$}
    \State{$\Delta = b^2 - 4 a c$}
    \If{$\Delta<0$}
    \State{\Return{$\emptyset$}}
    \EndIf
    \State{$Q = \lfloor \sqrt{\Delta}\rbrack$} \Comment{can be computed in integer arithmetic}
    \If{$(a<0) \vee (Q^2 \ne \Delta)$}
    \State{$I_0 = \lceil (-b-Q)/2 a\rceil$, $I_1 = \lfloor (-b+Q)/2 a\rfloor$}
    \Else
    \State{$I_0 = \lfloor (-b-Q)/2 a + 1\rfloor$,
      $I_1 = \lceil (-b+Q)/2 a - 1 \rceil$}
    \EndIf
    \If{$I_0\leq I_1$}
    \State{\Return{$\lbrack I_0, I_1\rbrack$}}
    \EndIf
    \State{\Return{$\emptyset$}}
    \EndFunction


    \vskip 3pt
    \noindent {\bf Time:} $O(1)$. \;\;\;\;\;\; {\bf Space:} $O(1)$.
        \end{algorithmic}
\end{algorithm}

\begin{algorithm}
  \caption{A very simple sieve of Eratosthenes}\label{alg:erasieve}
  \begin{algorithmic}[1]
    \Function{SimpleSiev}{$N$} 
    \Ensure{for $1\leq n\leq N$, $P_n=1$ if $n$ is prime, $P_n=0$ otherwise}
    \State{$P_1\gets 0$, $P_2\gets 1$, $P_n\gets 0$ for $n\geq 2$ even,
      $P_n\gets 1$ for $n\geq 3$ odd}
    \State{$m\gets 3$, $n\gets m\cdot m$}
    \While{$n\leq N$}
    \If{$P_m=1$}
    \While{$n\leq N$} \Comment{[sic]}
    \State{$P_n\gets 0$, $n\gets n + 2 m$} \Comment{sieves odd multiples
    $\geq m^2$ of $m$}
    \EndWhile
    \EndIf
    \State{$m\gets m+2$, $n\gets m\cdot m$}
    \EndWhile
    \State{\Return{$P$}}
    \EndFunction
    \vskip 3pt
    \noindent {\bf Time:} $O(N \log \log N)$.\; {\bf Space:} $O(N)$.
  \end{algorithmic}
\end{algorithm}

\begin{algorithm}
   \caption{A segmented sieve of Eratosthenes for finding primes}\label{alg:oldsegsieve}
  \begin{algorithmic}[1]
\Function{SegPrimes}{$n$,$\Delta$}
\Comment{finds all primes in $\lbrack n,n+\Delta\rbrack$}
\Ensure{ $S_j=\begin{cases} 1 &\text{if $n+j$ is prime}\\ 0 &\text{otherwise}\end{cases}$}
\State{$S_j \gets 1$ for all $0\leq j\leq \Delta$}
\State{$S_j\gets 0$ for $0\leq j\leq 1-n$} \Comment{[sic; excluding $0$ and $1$ from prime list]}
       \State{$M\gets \lfloor \sqrt{n+\Delta}\rfloor$, $P\gets \text{\textsc{SimpleSiev}}(M)$}
    \For{$1\leq m\leq M$}
    \If{$P_m=1$}
    \State{$n' \gets \max(m\cdot \lceil n/m\rceil,2 m)$}
    \While{$n' \leq n+\Delta$}\Comment{$n'$ goes over mults.\ of $m$ in $n+\lbrack 0,\Delta\rbrack$}
\State{$S_{n'-n} \gets 0$, $n' \gets n' + m$}
\EndWhile
\EndIf
\EndFor
\State{\Return{$S$}}
\EndFunction
\vskip 3pt
\noindent {\bf Time:} $O((\sqrt{n} + \Delta) \log \log (n+\Delta))$.\; {\bf Space:}
$O(n^{1/2} + \Delta)$.

\end{algorithmic}
\end{algorithm}

\begin{algorithm}
   \caption{A segmented sieve of Eratosthenes for computing $\mu(n)$}\label{alg:segsievemu}
  \begin{algorithmic}[1]
\Function{SegMu}{$n_0$,$\Delta$}
\Comment{computes $\mu(n)$ for $n$ in $\lbrack n_0,n_0+\Delta\rbrack$}
\Ensure{for $0\leq j\leq \Delta$,
  $m_j=\mu(n_0+j)$}
\State{$m_j \gets 1$, $\Pi_j\gets 1$ for all $0\leq j\leq \Delta$}
\State{$P\gets \text{\textsc{SimpleSiev}}(\lfloor\sqrt{n_0+\Delta}\rfloor)$}
\For{$p\leq \sqrt{n_0+\Delta}$}
\If{$P_p=1$} \Comment{if $p$ is a prime\dots}
\State{$n \gets p\cdot \lceil n_0/p\rceil$}\Comment{smallest multiple $\geq n_0$ of $p$}
\While{$n \leq n_0+\Delta$} \Comment{$n$ goes over multiples of $p$}
\State{$m_{n-n_0} \gets -m_{n-n_0}$, $\Pi_{n-n_0} = p\cdot \Pi_{n-n_0}$, $n \gets n + p$}
\EndWhile
\State{$n \gets p^2\cdot \lceil n_0/p^2\rceil$}\Comment{smallest multiple $\geq n_0$ of $p^2$}
\While{$n \leq n_0+\Delta$} \Comment{$n$ goes over multiples of $p^2$}
\State{$m_{n-n_0} \gets 0$, $n \gets n + p^2$}
\EndWhile
\EndIf
\EndFor
\For{$0\leq j\leq \Delta$}
\If{$m_j\ne 0 \wedge \Pi_j\ne n_0+j$}
\State{$m_j \gets -m_j$}
\EndIf
\EndFor
\State{\Return{$m$}}
\EndFunction
\vskip 3pt
\noindent {\bf Time:} $O((\sqrt{n_0} + \Delta) \log \log (n_0+\Delta))$.\; {\bf Space:}
$O(\sqrt{n_0} + \Delta \log (n_0+\Delta))$,
or, after a standard improvement (\S \ref{sec:impdet}), space
$O(\sqrt{n_0} + \Delta \log \log (n_0+\Delta))$.
  \end{algorithmic}
  \end{algorithm}

\begin{algorithm}
  \caption{A segmented sieve of Eratosthenes for factorization}\label{alg:basterno}
  \begin{algorithmic}[1]
\Function{SubSegSievFac}{$n$,$\Delta$,$M$}
\Comment{finds prime factors $p\leq M$}
\Ensure{for $0\leq j\leq \Delta$, $F_j=\{(p,v_p(n+j))\}_{p\leq M, p|n+j}$}
\Ensure{for $0\leq j\leq \Delta$, $\Pi_j = \prod_{p\leq M, p|(n+j)} p^{v_p(n+j)}$.}
\State{$F_j \gets\emptyset$, $\Pi_j\gets 1$ for all $0\leq j\leq \Delta$}
\State{$\Delta'\gets \lfloor \sqrt{M}\rfloor$, $M'\gets 1$}
\While{$M'\leq M$}
\State{$P\gets \text{\textsc{SegPrimes}}(M',\Delta')$}
\For{$M'\leq p<M'+\Delta'$}
\If{$P_{p-M'}=1$}  \Comment{if $p$ is a prime\dots}
\State{$k\gets 1$, $d\gets p$} \Comment{$d$ will go over the powers $p^k$ of $p$}
\While{$d\leq n+\Delta$}
\State{$n'\gets d\cdot \lceil n/d\rceil$}
\While{$n'<x$}
\If{$k=1$}
\State{{\bf append} $(p,1)$ to $F_{n'-n}$}
\Else
\State{{\bf replace} $(p,k-1)$ by $(p,k)$ in $F_{n'-n}$}
\EndIf
\State{$\Pi_{n'-n}\gets p\cdot \Pi_{n'-n}$, $n'\gets n'+d$}
\EndWhile
\State{$k\gets k+1$, $d\gets p\cdot d$}
\EndWhile
\EndIf
\EndFor
\State{$M'\gets M'+\Delta'$}
\EndWhile
\State{\Return{$(F,\Pi)$}}
\EndFunction
\vskip 3pt
\noindent {\bf Time:} $O((M + \Delta) \log \log (n+\Delta))$,
\vskip 1pt
\noindent {\bf Space:}
$O(M + \Delta \log (n+\Delta))$.
 \end{algorithmic}
\end{algorithm}

\begin{algorithm}
  \caption{A segmented sieve of Eratosthenes for factorization, II}\label{alg:basterno2}
\begin{algorithmic}[1]
\Function{SegFactor}{$n$,$\Delta$}
\Comment{factorizes all $n'\in \lbrack n, n+\Delta\rbrack$}
\Ensure{for $0\leq j\leq \Delta$, $F_j$ is the list of pairs $(p,v_p(n+j))$
  for $p|n+j$}
\State{$(F,\Pi) \gets \textsc{SubSegSievFac}(n,\Delta,\lfloor \sqrt{x}
  \rfloor)$}
\For{$n\leq n'\leq n+\Delta$}
\If{$\Pi_{n'-n} \ne n'$}
\State{$p_0\gets n'/\Pi_{n'-n}$, {\bf append} $(p_0,1)$ to $F_{n'-n}$}
\EndIf
\EndFor
\State{\Return{$F$}}
\EndFunction
\vskip 3pt
\noindent {\bf Time:} $O((\sqrt{n} + \Delta) \log \log (n+\Delta))$,
{\bf Space:}
$O(\sqrt{n} + \Delta \log (n+\Delta))$.
 \end{algorithmic}
\end{algorithm}

\begin{algorithm}
  \caption{From factorizations to $\sum_{d|n: d\leq a} \mu(d)$}\label{alg:factofun}
\begin{algorithmic}[1]

\vskip 5pt
\Function{SubFacTSM}{$F$, $m$,$m'$,$a$,$n$}
\If{$m>a$}
\State{\Return{$0$}}
\EndIf
\If{$F=\emptyset$}
\State{\Return{$1$}}
\EndIf
\If{$m' a \geq n$}
\State{\Return{$0$}}
\EndIf
\State{\textbf{Choose} $(p,i)\in F$ such that $p$ is maximal}
\State{$F' = F\setminus \{(p,i)\}$}
\State{\Return{\textsc{SubFacTSM}($F'$,$m$,$p m'$,$a$,$n$) -
    \textsc{SubFacTSM}($F'$,$m p$,$m'$,$a$,$n$)}}
\EndFunction
    
\vskip 5pt
\Function{FacToSumMu}{$F$, $a$}
\Require{$F$ is the list of all pairs $(p,v_p(n))$, $p|n$, for some $n$, with $p$ in order}
\Ensure{returns $\sum_{d|n: d\leq a} \mu(d)$}
\State{$n' = \prod_{(p,i)\in F} p$}
\State{\Return{\textsc{SubFacTSM($F$,$1$,$1$,$a$,$n'$)}}}
\EndFunction
\vskip 3pt
\noindent {\bf Time:} $O(2^{\len(F)})$, but less on average (see Prop \ref{prop:askesis}). {\bf Space:} $O(\len(F))$.
  \end{algorithmic}
\end{algorithm}

 

\clearpage

\bibliographystyle{alpha}
\bibliography{mertensformic}
\end{document}